\documentclass{amsart}
\usepackage[OT1]{fontenc}
\usepackage{amsthm,amsmath,amsfonts,amssymb}
\usepackage[all,cmtip]{xy}
\usepackage[pdftex]{graphicx}
\usepackage[colorlinks]{hyperref}
\usepackage{tikz}
\usetikzlibrary{positioning}

\newtheorem{thm}{Theorem}[section]
\newtheorem{lem}[thm]{Lemma}
\newtheorem{prop}[thm]{Proposition}
\newtheorem{cor}[thm]{Corollary}

\theoremstyle{definition}
\newtheorem{defn}[thm]{Definition}

\newtheorem{exmp}[thm]{Example}

\theoremstyle{remark}
\newtheorem{rem}[thm]{Remark}

\newcommand{\cA}{\mathcal{A}}

\newcommand{\C}{\mathbb{C}}
\newcommand{\cC}{\mathcal{C}}

\newcommand{\E}{\mathbb{E}}

\newcommand{\cH}{\mathcal{H}}

\newcommand{\cI}{\mathcal{I}}

\newcommand{\cL}{\mathcal{L}}
\newcommand{\cM}{\mathcal{M}}
\newcommand{\mm}{\mathfrak{m}}

\renewcommand{\P}{\mathbb{P}}
\newcommand{\cP}{\mathcal{P}}

\newcommand{\R}{\mathbb{R}}

\newcommand{\cT}{\mathcal{T}}

\newcommand{\cX}{\mathcal{X}}

\newcommand{\indep}{{\;\bot\!\!\!\!\!\!\bot\;}}
\DeclareMathOperator{\indic}{{1\!\!1}}

\begin{document}

\title{$L$-cumulants, $L$-cumulant embeddings and algebraic statistics}
\author{{Piotr} {Zwiernik}}
\markboth{P. Zwiernik}{$L$-cumulants and product partition models}
\address{Piotr Zwiernik\\TU Eindhoven\\ Department of Mathematics and Computer Science\\ PO Box 513\\ 5600 MB
Eindhoven\\ The Netherlands}
\email{piotr.zwiernik@gmail.com}
\thanks{}
\date{}

\begin{abstract} Focusing on the discrete probabilistic setting we generalize the combinatorial definition of cumulants to $L$-cumulants. This generalization keeps all the desired properties of the classical cumulants like semi-invariance and vanishing for independent blocks of random variables. These properties make $L$-cumulants useful for the algebraic analysis of statistical models. We illustrate this for general Markov models and hidden Markov processes in the case when the hidden process is binary. The main motivation of this work is to understand cumulant-like coordinates in algebraic statistics and to give a more insightful explanation why tree cumulants give such an elegant description of binary hidden tree models. Moreover, we argue that $L$-cumulants can be used in the analysis of certain classical algebraic varieties. \end{abstract}
\maketitle
\section{Introduction}\label{sec:introduction}

Although moments provide a convenient summary of properties of a probability distribution, it was observed that these properties can generally be  described in a simpler way using cumulants (see for example \cite[Section 2.4]{davison2003}, \cite[Chapter 2]{mccullagh1987tms}). This is mainly because cumulants have the ability to capture symmetries and underlying independencies of a probability distribution. These striking features of cumulants make them an interesting object of statistical study both from a theoretical and practical point of view. In addition, as it was shown for example in \cite{geiger2001sef,rusakov2006ams,settimi2000gma}, cumulants and moments can be used to analyze the geometry of statistical models. 

Recently, in  \cite{pwz-2010-identifiability} we have suggested using a less standard system of coordinates which we called \textit{tree cumulants}. This new coordinate system proved to be useful to analyze Bayesian networks on trees when some of the nodes are not observed. Various results on identifiability and geometry of these models have been obtained in \cite{pwz-2010-bic,pwz-2009-semialgebraictrees,pwz-2010-identifiability}, which encouraged us to study more general coordinate systems like that. In the present paper we propose a useful generalization of both cumulants and tree cumulants. 

We work in a simple probabilistic setting. Let $X=(X_{1},\ldots,X_{n})$ be a random vector such that each $X_{i}$ takes $r_{i}\geq 2$ possible values, where each $r_i$ is finite. The vector $X$ takes values in a finite discrete set $\cX=\prod_{i=1}^n\cX_i\subseteq \R^n$ such that $|\cX_i|=r_{i}$ for $i=1,\ldots,n$. Without loss of generality we set
$$\cX\quad=\quad\{0,\ldots,r_{1}-1\}\,\,\times\,\, \cdots\,\,\times\,\, \{0,\ldots, r_{n}-1\}.$$
Any probability distribution of $X$ can be written as a point $P=[p(x)]\in \mathbb{R}^\cX$ such that $p(x)\geq 0$ for all $x\in \cX$ and $\sum_{x\in \cX} p(x)=1$. The set of all such points is called the \textit{probability simplex} and it is denoted by $\Delta_\cX$. 

For any function $f:\, \cX\rightarrow \R$ the \textit{expectation of $f(X)$} is given by
$$
\E[f(X)]\quad:=\quad\sum_{x\in \cX}p(x)f(x).
$$
Let $[n]:=\{1,\ldots,n\}$ and for any multiset $A=\{i_1,\ldots,i_d\}$ of elements of $[n]$ let $$X_A=(X_{i_1},\ldots,X_{i_d}).$$ In a similar way we define $x_A=(x_{i_1},\ldots,x_{i_d})$ and $\cX_A=\cX_{i_1}\times \cdots\times \cX_{i_d}$. For each such a multiset $A$ we define the corresponding \textit{moment}
$$
\mu_A\quad=\quad\E[X_{i_1}\cdots X_{i_d}]
$$ 
and the \textit{central moment}
$$
\mu'_A\quad=\quad\E[(X_{i_1}-\mu_{i_1})\cdots (X_{i_d}-\mu_{i_d})].
$$ 
Our convention is to write $\mu_A$ as $\mu_{i_1\cdots i_d}$, where $i_1\leq \cdots\leq i_d$. So for example if $A=\{1,2,4,4,4\}$, the corresponding moment is written as $\mu_{12444}=\E[X_1X_2X_4^3]$. The same convention applies to central moments. In particular, for every $i<j$, $\mu'_{ij}$ is the covariance between $X_i$ and $X_j$.

To show how cumulants can be naturally generalized we first define them formally and then we discuss their basic properties. Cumulants are usually computed using the cumulant generating function, which is defined as the logarithm of the moment generating function.  In this paper we use an alternative definition of cumulants using partitions (see for example \cite{mccullagh1987tms,rota_cumulants,speed1983cumulants}). We say that $\pi=B_{1}|\ldots|B_{k}$ is a \textit{partition} (or a \textit{set partition}) of $[n]$, if the \textit{blocks} $B_{i}\neq \emptyset$ are  disjoint sets whose union is $[n]$. A partition is called a \textit{split} if it consists of two blocks. Let $\Pi([n])$ be the set of all set partitions of $[n]$.  The \textit{cumulant of the vector $X$} is defined as
\begin{equation}\label{eq:cum1}
{k}_{1\cdots n}\quad=\quad\sum_{\pi\in\Pi([n])}(-1)^{|\pi|-1}(|\pi|-1)!\prod_{B\in \pi} \mu_{B},
\end{equation}
where the sum is over all set partitions of $[n]$, the product is over all blocks of a partition and $|\pi|$ denotes the number of blocks of $\pi$. For example, if $n=3$ then there are five partitions in $\Pi({[3]})$: $123$, $1|23$, $2|13$, $12|3$ and $1|2|3$ and (\ref{eq:cum1}) gives
\begin{equation}\label{eq:kappa3}
{k}_{123}\quad=\quad\mu_{123}-\mu_{1}\mu_{23}-\mu_{2}\mu_{13}-\mu_{12}\mu_{3}+2\mu_{1}\mu_{2}\mu_{3}.
\end{equation}
Equation (\ref{eq:cum1}) can be generalized for any multiset $A=\{i_{1},\ldots, i_{d}\}$ of elements of $[n]$ to obtain the cumulant of $X_A$. We use the bijection between $A$ and $[d]$ and write 
\begin{equation}\label{eq:cumA}
{k}_{A}\quad=\quad\sum_{\pi\in\Pi([d])}(-1)^{|\pi|-1}(|\pi|-1)!\prod_{B\in \pi} \mu_{i_B},
\end{equation}
where $i_B=\{i_j:\,j\in B\}$. Hence for instance
$$
{k}_{112}\quad=\quad\mu_{112}-2\mu_{1}\mu_{12}-\mu_{11}\mu_{2}+2\mu_{1}^{2}\mu_{2}.
$$

For each $x=(x_1,\ldots,x_n)\in \cX$ define a multiset $\cA(x)$ as 
\begin{equation}\label{eq:defAX}
\cA(x)\quad=\quad\{\underbrace{1,\ldots,1}_{{x_{1}\,{\tiny\mbox{times}}}},\ldots, \underbrace{n,\ldots,n}_{{x_{n}\,{\tiny\mbox{times}}}}\}
\end{equation}
and let $\cA(\cX)=\{\cA(x):\,x\in \cX\}$. By the \textit{moment aliasing principle} (see \cite[Lemma~3]{pistone2006cv}) there exists a polynomial isomorphism between $P=[p(x)]_{x\in \cX}$ and two other systems of coordinates of $\R^{\cA(\cX)}\simeq \R^\cX$ given by moments $M=[\mu_{\cA(x)}]_{x\in \cX}$ and by cumulants $K=[k_{\cA(x)}]_{x\in \cX}$. In particular, every model $\cM\subseteq \Delta_{\cX}$, after a change of coordinates, can be equivalently expressed in terms of $M$ or $K$. 

In our discussion of cumulants the central concept is that of independence. Let $B\subseteq [n]$ and define the function $\indic_{x_B}$ on $\cX$ by  $\indic_{x_B}(X)=1$ if $X_B=x_B$ and $\indic_{x_B}(X)=0$ otherwise. By $p_B$ denote the \textit{marginal distribution of $X_B$} defined by 
$$
p_B(x_B)\quad=\quad\E[\indic_{x_B}(X)]\qquad\mbox{for every }x_B\in \cX_B.
$$
For any two disjoint subsets $I,J\subseteq [n]$ we say that \textit{$X_I$ and $X_J$ are independent}, which we denote by $I\indep J$ (or $X_I\indep X_J$), if and only if 
$$p_{I\cup J}(x_{I\cup J})=p_I(x_I)p_J(x_J)\qquad \mbox{for all } x\in \cX.$$ 
The following formulation of independence in terms of moments will be helpful. 
\begin{lem}\label{lem:indinmom}
We have $I\indep J$ for some disjoint sets $I,J\subseteq[n]$ if and only if 
$$
\mu_{A\cup B}\,\,=\,\,\mu_{A}\mu_{B}\qquad\mbox{for all nonempty } A\in \cA(\cX_I), B\in \cA(\cX_J),
$$
where $\cA(\cX_I)=\{\cA(x):\,x\in \cX_I\}$.
\end{lem}
\begin{proof}We use an alternative definition of independence (see \cite[page 136]{feller1971ipt}) which states that $X_I$ and $X_J$ are independent if and only if for any two $L^2$-functions $f,g$ we have $$\E[f(X_I)g(X_J)]=\E[f(X_I)]\E[g(X_J)].$$ Now the 'if' direction of the lemma is immediate. The 'only if' direction uses the fact that the set of values of $X$ is discrete and finite. In this case any function of $X$ is a polynomial function (can be represented as a polynomial in the entries of $X$), where the terms of these polynomials are $\prod_{i\in A}X_i$ for all $A\in \cA(\cX)$. Thus, to check if $I\indep J$, it remains to check if
$$
\E [f(X_{I})g(X_{J})]\quad=\quad\E [f(X_{I})]\E [g(X_{J})]
$$  
for all polynomials $f,g$ such that each $f$ has only terms $\prod_{i\in A}X_i$ for all nonempty \mbox{$A\in \cA(\cX_{I})$} and $g$ has only terms $\prod_{i\in B}X_i$ for all $B\in \cA(\cX_{J})$. By expanding the terms of $f$ and $g$ it suffices to check that this property holds for each monomial, which is true by assumption.
\hfill$\Box$\end{proof}

\begin{exmp}Let $m=2$, $r_1=2$ and $r_2=3$. Then $\cX=\{0,1\}\times \{0,1,2\}$ and $$\cA(\cX)=\{\emptyset,\{2\},\{2,2\},\{1\},\{1,2\},\{1,2,2\}\}.$$
Since $\cA(\cX_1)=\{\emptyset,\{1\}\}$ and $\cA(\cX_2)=\{\emptyset,\{2\},\{2,2\}\}$, by Lemma \ref{lem:indinmom}, we have $1\indep 2$ if and only if $\mu_{12}=\mu_1\mu_2$, $\mu_{122}=\mu_1\mu_{22}$, where $\mu_{122}=\E[X_1X_2^2]$.
\end{exmp}

Lemma \ref{lem:indinmom} generalizes and we have  $I_1\indep I_2\indep \cdots \indep I_r$ for some disjoint sets $I_1,\ldots, I_r\subseteq[n]$ if and only if
\begin{equation}\label{eq:indepinmoms}
\mu_{A_1\cdots A_r}\,\,\,=\,\,\,\prod_{i=1}^r \mu_{A_i},\qquad \mbox{for all }A_i\in \cA(\cX_{I_i}),\, i=1,\ldots,r,
\end{equation}
where $A_1\cdots A_r$ is a shorthand notation for $A_1\cup \cdots \cup A_r$.

Cumulants satisfy the following four basic properties, which make them useful for statistical modelling.
\begin{itemize}
\item[(P1)] Whenever there exists a split of the set of indices $[n]$ of $X$ into two block $A|B$ such that  $A\indep B$ then ${k}_{1\cdots n}=0$. 
\item[(P2)] For any $a\in \R^{n}$ define $\widetilde{X}=X+a$ and for any multiset $A$ by $\widetilde{k}_A$ denote the corresponding cumulant of $\widetilde{X}_A$. Then $\widetilde{k}_i=k_i+a_{i}$ for every $i=1,\ldots,n$, and $\widetilde{k}_A=k_A$ whenever $|A|\geq 2$. 
\item[(P3)] Let $Q=[q_{ij}]\in\R^{m\times n}$, $X\in \R^n$ and let $\widetilde{X}=QX\in \R^m$. Define $\widetilde{k}_A$ as the cumulant of $\widetilde{X}_A$, where $A$ is a multiset of elements of $[m]$. Let $K^{(d)}=[k_{i_1\cdots i_d}]$ denote the $(n\times \cdots \times n)$-tensor indexed by all multisets of elements of $[n]$ of size $d\geq 1$; and let $\widetilde{K}^{(d)}=[\widetilde{k}_{i_1\cdots i_d}]$ be the $(m\times \cdots \times m)$-tensor indexed by all multisets of elements of $[m]$. Then, 
$$\widetilde{K}^{(d)}\quad=\quad Q\cdot K^{(d)},\qquad\mbox{for every }d\geq 1$$ where for every multiset $\{i_1,\ldots,i_d\}$ of elements of $[m]$:
$$
(Q\cdot K^{(d)})_{i_{1}\cdots i_{d}}\quad:=\quad\sum_{j_{1}=1}^{n}\cdots\sum_{j_{d}=1}^{n}q_{i_{1}j_{1}}\cdots q_{i_{d}j_{d}} k_{j_{1}\cdots j_{d}}.
$$
In other words cumulants under linear mappings transform as contravariant tensors.
\item[(P4)] For two random vectors $X$, $Y$ of dimension $n$ denote by ${k}_A(X)$, ${k}_A(Y)$ and ${k}_A(X+Y)$ the cumulants of $X$, $Y$ and $X+Y$ respectively. If $X\indep Y$ then $k_A(X+Y)=k_A(X)+k_A(Y)$ for every multiset $A$ of elements of $[n]$.
\end{itemize}

In this paper we generalize cumulants by changing the set $\Pi([n])$ in (\ref{eq:cum1}) for other set partition lattices. The term $(-1)^{|\pi|}(|\pi|-1)!$ in each summand of (\ref{eq:cum1}) is replaced by another function of $\pi$ which will be specified later. These generalized cumulants keep usually all properties (P1)-(P4) of classical cumulants. Also the Brillinger's conditional cumulants formula derived in  \cite{brillinger1969calculation} can be generalized under additional conditions. 

Different forms of cumulants are known to researchers in  non-commutative probability. For example \textit{free cumulants} are used in the theory of random matrices \cite{lehner2004cumulants,speicher_freeprob97} and  \textit{Boolean} \textit{cumulants} are applied to stochastic differential equations \cite{Maruyama1988447}. All those cumulants fall under our general definition. In Proposition \ref{prop:centralm} we show that central moments can be also represented as generalized cumulants. As an interesting implication we get a simple computationally efficient formula for central moments in terms of moments (see Lemma \ref{lem:centrmomform}). The proof of this formula is straightforward. 

As it has been already pointed out in \cite{pwz-2011-bincum}, cumulants and cumulant-like quantities are also useful in algebraic geometry. The  coordinate system given by cumulants has a number of useful properties. For example, the tangential variety ${\rm Tan}((\P^1)^n)$, when expressed in binary cumulants, becomes toric. Also, the study of the secant variety ${\rm Sec}((\P^1)^n)$ becomes easier when we change coordinates to binary tree cumulants. This happens because the induced parametrization in this new coordinate system becomes nearly monomial (see Section~\ref{sec:tressforsecs}).

There are two main reasons why cumulants can be successfully used in algebraic geometry and in the geometric study in statistics. First, many interesting algebraic varieties coincide with some statistical models. Second, the whole machinery of cumulants is purely algebraic in the sense that nonnegativity of probabilities does not play any role. In fact the only condition which we impose on probabilities is that they sum to one. For that reason the same techniques can be applied to any complex tensor with coordinates summing to one. This observation links our work to the theory of umbral calculus  \cite{rota1973foundations}. 

This paper is organized as follows. In Section \ref{sec:XClattice} we introduce some basic concepts of the theory of partially ordered sets. In Section \ref{sec:binarycums} we define binary $L$-cumulants, which form a rather straightforward generalization of binary cumulants introduced in \cite{pwz-2011-bincum}. In Section \ref{sec:geometry} we present how binary $L$-cumulants may be used in algebraic geometry. This is then exemplified with a basic study of secant varieties in Section \ref{sec:tressforsecs}. The general definition of $L$-cumulants is provided in Section \ref{sec:XCdef}. In Section \ref{sec:XCproperties} we show that, under some mild conditions,  all the basic properties (P1)-(P4) of classical cumulants hold also for $L$-cumulants. Moreover, in Section \ref{sec:XCconditional} we generalize the Brillinger's formula for cumulants in terms of conditional cumulants. In Section \ref{sec:treecum} we show how the results of this paper explain why tree cumulants work so well for tree models. We also provide a simple analysis of processes with an underlying hidden two-state Markov chain, which in particular gives a very simple parametrization of homogeneous binary hidden Markov models.

\section{Basic combinatorics}\label{sec:XClattice}

In this section we introduce basic combinatorial concepts used later in the paper. For a more detailed treatment see \cite{stanley2006enumerative}. Recall that $\pi=B_{1}|\ldots|B_{k}$ is a partition of $[n]$, if the \textit{blocks} $B_{i}\neq \emptyset$ are  disjoint sets whose union is $[n]$. Equivalently, a partition of $[n]$ corresponds to an equivalence relation $\sim_{\pi}$ on $[n]$ where $i\sim_{\pi} j$ if $i$ and $j$ lie in the same block.    Let now $A$ be a multiset $A=\{i_{1},\ldots,i_{d}\}$ of elements of $[n]$. We define a partition $\pi$ of $A$ using a partition $\pi$ of $[d]$ by $i_{j}\sim_{\pi} i_{k}$ if $j\sim_{\pi} k$ in $\Pi([d])$. The set of all partitions of $A$ is denoted by $\Pi(A)$ and by definition it is isomorphic to $\Pi([d])$. 

A \textit{partially ordered set $\cP$} (or \textit{poset}) is a set together with an ordering $\leq$ such that: $\pi\leq \pi$ for all $\pi\in \cP$; if $\pi\leq \nu$ and $\nu\leq \pi$ then $\pi=\nu$; and if $\pi\leq \nu$ and $\nu\leq \delta$ then $\pi\leq \delta$ for all $\pi,\nu,\delta\in \cP$. A \textit{subposet} of $\cP$ is any subset of $\cP$ with the same ordering.  As an important example of a poset consider the set $\Pi([n])$ with the poset structure given by refinement ordering such that $\pi\leq \nu$ in $\Pi([n])$ if and only if every block of $\pi$ is contained in a block of $\nu$. For instance let $n=5$, $\pi=13|4|25$ and $\nu=1235|4$ then $\pi\leq \nu$. 

We say that \textit{$\cP$ has a $\hat{0}$} if there exists an element $\hat{0}\in \cP$ such that $\pi\geq \hat{0}$ for all $\pi\in \cP$. Similarly, \textit{$\cP$ has a $\hat{1}$} if there exists $\hat{1}\in \cP$ such that $\pi\leq \hat{1}$ for all $\pi\in \cP$. If $\pi$ and $\nu$ belong to a poset $\cP$, then an \textit{upper bound} of $\pi$ and $\nu$ is an element $\delta\in \cP$ satisfying $\delta\geq \pi$ and $\delta\geq \nu$. A \textit{least upper bound} of $\pi$ and $\nu$ is an upper bound $\delta$ of $\pi$ and $\nu$ such that every upper bound $\gamma$ of $\pi$ and $\nu$ satisfies $\gamma\geq \delta$. If a least upper bound of $\pi$ and $\nu$ exists, then it is clearly unique and it is denoted by $\pi\vee \nu$. Dually one can define the \textit{greatest lower bound} $\pi\wedge \nu$ when it exists. We call $\vee$ the \textit{join operator} and $\wedge$ the \textit{meet operator}. 

A \textit{lattice} is a poset $L$ for which every pair of elements has a least upper bound and greatest lower bound. A \textit{sublattice} of a lattice $L$ is a nonempty subset of $L$ which is a lattice with \textit{the same} meet and join operations as $L$. Clearly all finite lattices have a $\hat{0}$ and $\hat{1}$. In particular $\Pi([n])$ forms a lattice where the $n$-block partition $1|2|\cdots|n$ is the $\hat{0}$, and the one-block partition $[n]$ is the $\hat{1}$ of this lattice. A \textit{meet semilattice} is a poset $S$ for which every pair of elements has a least upper bound. A \textit{meet subsemilattice of $S$} is a subposet of $S$ which forms a meet semilattice with the same meet operator as $S$. Dually we define a \textit{join semilattice} and a \textit{join subsemilattice}.

\begin{defn}By a \textit{partition lattice} of a set $[n]$ we mean any lattice $L$ which forms a subposet of $\Pi([n])$ and both the one block partition $[n]$ and the minimal partition $1|2|\cdots|n$ lie in $L$. 
\end{defn}

Note that we do not require that a partition lattice forms a sublattice of $\Pi([n])$. 

\begin{defn}\label{def:partitions} The following is a list of interesting set partition lattices.
\begin{itemize}
\item[(1)] A partition $\pi\in \Pi([n])$ is \textit{non-crossing} if there is no quadruple of elements $i<j<k<l$ such that $i\sim_{\pi} k$, $j\sim_{\pi} l$ and $i\nsim_{\pi}j$. The {noncrossing partitions of $[n]$} form a lattice which we denote by ${\rm NC}([n])$. This lattice is not a sublattice of $\Pi([n])$, however, it is a meet subsemilattice of $\Pi([n])$ because the meet operators coincide.
\item[(2)] An \textit{interval partition} of $[n]$ is a partition $\pi$ of a form $$1\cdots i_{1}|(i_{1}+1)\cdots i_{2}|\cdots|(i_{k}+1)\cdots n$$ for some $0\leq k\leq n-1$ and $1\leq i_{1}<\ldots<i_{k}\leq n-1$. The poset of all interval partitions is denoted by $\cI([n])$. It forms a sublattice of $\Pi([n])$  isomorphic to the Boolean lattice of $[n-1]$.
\item[(3)] A partition $\pi\in \Pi([n])$ is called a \textit{one-cluster partition} if it contains at most one block of size greater than one. In particular the one-block partition $[n]$ and the minimal partition $1|2|\cdots|n$ are one-cluster partitions. The poset of all one-cluster partitions forms a lattice $\cC([n])$, which is not a sublattice of $\Pi([n])$. It is isomorphic to the poset of all subsets of $[n]$ excluding singletons.  It forms a meet subsemilattice of $\Pi([n])$. 
\item[(4)] Let $T=(V,E)$ be a fixed tree with set of nodes $V$, set of edges $E$ and with $n$ leaves labelled by $[n]$. Removing a subset of edges $E'$ from $E$ induces a forest. Restricting $[n]$ to the connected components of this forest gives a \textit{tree partition} $\pi$ induced by $T$. The set of all tree {partitions} induced by $T$ is denoted by $\cT^{T}([n])$ and it forms a lattice which is a meet subsemilattice of $\Pi([n])$. For an example of a tree and the induced lattice of partitions see Figure \ref{fig:cater} (for $n=4$) and  Figure \ref{fig:treeposet4}.
\end{itemize}
\end{defn}

For every poset $\cP$ we define the \textit{M\"{o}bius function $\mm_{\cP}:\cP\times\cP\rightarrow \R$} by 
\begin{equation}\label{eq:mobius-def}
\mm_{\cP}(\pi,\nu)\,\,=\,\,\left\{\begin{array}{ll}
 1,& \mbox{if } \pi=\nu,\\
-\sum_{\pi\leq \delta<\nu} \mm_{\cP}(\pi,\delta) &\mbox{if } \pi<\nu, \\
 0,&\mbox{otherwise.}
\end{array}\right.
\end{equation}
When there is no ambiguity we usually drop $\cP$ in the notation denoting the M\"{o}bius function on $\cP$ by $\mm$. Note that directly from the definition in (\ref{eq:mobius-def})
\begin{equation}\label{eq:sumintzero}
\sum_{\pi\leq\delta\leq \nu} \mm_{\cP}(\pi,\delta)\,\,\,=\,\,\,\left\{\begin{array}{l}
0\quad\mbox{if }\pi<\nu\\
1\quad\mbox{if }\pi=\nu.
\end{array}\right.
\end{equation}

A special type of a subposet of $\cP$ is the \textit{interval} $$[\pi,\nu]=\{\delta\in\cP:\, \pi\leq \delta\leq \nu\},$$ defined whenever $\pi\leq \nu$. The M\"{o}bius function on this subposet is naturally induced from the  M\"{o}bius function on $\cP$ (see for example \cite[Proposition~4]{rota1964fct}). For any two posets $\cP_{1},\cP_{2}$ we define the poset $\cP_{1}\times\cP_{2}$ as a set with the ordering $(\pi,\nu)\leq (\pi',\nu')$ if $\pi\leq \pi'$ and $\nu\leq\nu'$. The following result gives a convenient way of finding a M\"{o}bius function for posets constructed from other posets by taking products. 
\begin{prop}[Proposition~3.1.2,  \cite{stanley2006enumerative}]\label{prop:mobius-product} Let $\cP_{1}$ and $\cP_{2}$ be  finite posets, and let $\cP_{1}\times\cP_{2}$ be their direct product. If $(\pi,\nu)\leq (\pi',\nu')$ in $\cP_{1}\times\cP_{2}$, then
$$
\mm_{\cP_{1}\times\cP_{2}}((\pi,\nu),(\pi',\nu'))\quad=\quad\mm_{\cP_{1}}(\pi,\nu)\mm_{\cP_{2}}(\pi',\nu').
$$
\end{prop}
The M\"{o}bius function is especially useful due to the following result.
\begin{prop}[M\"{o}bius inversion formula]\label{prop:mobdual} Let $\cP$ be a finite poset. Let $f,g:\cP\rightarrow\R$. Then 
$$
g(\pi)\quad=\quad\sum_{\nu\leq \pi} f(\nu),\quad\mbox{for all } \pi\in \cP,
$$
if and only if 
$$
f(\pi)\quad=\quad\sum_{\nu\leq \pi}\mm(\nu,\pi)\,g(\nu)\quad\mbox{for all } \pi\in \cP.
$$
\end{prop}
For every lattice denote $\mm(\pi):=\mm(\pi,\hat{1})$. Later we will see that it is particularly important to identify values of $\mm(\pi)$ for various partition lattices. For $\Pi([n])$ we have $\mm(\pi)=(-1)^{|\pi|-1}(|\pi|-1)!$ The lattice of interval partitions $\cI([n])$ is isomorphic to the Boolean lattice of all subsets of $[n-1]$ and hence $\mm(\pi)=(-1)^{|\pi|-1}$. For the lattice of one-cluster partitions we have
\begin{equation}\label{eq:Mobiusfo1c}
\mm(\pi)\quad=\quad\left\{\begin{array}{ll}
(-1)^{n-1}(n-1) & \mbox{if } \pi=1|2|\cdots|n, \mbox{ and}\\
(-1)^{|\pi|-1} & \mbox{otherwise}.
\end{array}
\right.
\end{equation}
For the other cases in Definition~\ref{def:partitions} the M\"{o}bius function can be computed recursively.

\section{Binary $L$-cumulants}\label{sec:binarycums}

In this section we discuss binary $L$-cumulants which generalize binary cumulants of \cite{pwz-2011-bincum}. Most of the technical results will be stated without proofs, which will then be given in a more general context in later sections.

\subsection{Definition~and basic facts}
Assume that $\cX=\{0,1\}^n$, in which case $\cA(\cX)$ is the set of all subsets of $[n]$. Let $L\subseteq\Pi([n])$ be a partition lattice of $[n]$. For every $I\subseteq [n]$ consider $L(I)$ as the subposet of $\Pi(I)$ obtained from $L$ by constraining each partition to the subset $I$. The M\"{o}bius function on $L(I)$ is also denoted by $\mm$ unless it may lead to ambiguity in which case we write explicitly $\mm_I$.

 A \textit{multiplicative function} on $L(I)$ is any function such that for every $\pi\in L(I)$ 
$$
f(\pi)\quad=\quad\prod_{B\in \pi} f_{B}\qquad\mbox{for some } f_{B}\in \R.
$$
First consider the case when $L=\Pi([n])$. For every $I\subseteq[n]$ and $\nu\in \Pi(I)$ define
\begin{equation}\label{eq:cum2}
k(\nu)\quad=\quad \sum_{\pi\leq\nu}\mm(\pi,\nu)\mu(\pi),
\end{equation}
where $\mu(\pi)=\prod_{B\in \pi}\mu_{B}$ is a multiplicative function and the sum is taken over elements $\pi$ of $\Pi(I)$ such that $\pi\leq \nu$.  The one-block partition $I$ is the unique maximal element of the lattice $\Pi(I)$. The M\"{o}bius function on $\Pi(I)$ satisfies $\mm(\pi):=\mm(\pi,I)=(-1)^{|\pi|-{1}}(|\pi|-1)!$ for all $\pi\in \Pi(I)$. It follows by (\ref{eq:cumA}) that $k_I=k(I)$ and hence (\ref{eq:cum2}) evaluated at $\nu={I}$ gives the definition of binary cumulants. 

To get the inverse formula for moments in terms of cumulants we need the following result.
\begin{lem}\label{lem:kisprod}For every $\nu\in \Pi(I)$ we have $k(\nu)=\prod_{B\in \nu}k_B$, where $k(\nu)$ is defined by (\ref{eq:cum2}).
\end{lem}
\begin{proof}
Note that every interval $[\pi,\nu]\subseteq\Pi(I)$ is isomorphic to a product of intervals $\prod_{B\in \nu}[\pi(B),B]\subseteq \prod_{B\in \nu}\Pi(B)$, where $\pi(B)$ denotes $\pi$ constrained to elements in $B\subseteq I$. By Proposition~\ref{prop:mobius-product} a M\"{o}bius function on a product of posets is equal to the product of M\"{o}bius functions for each individual factor. Hence, (\ref{eq:cum2}) can be rewritten as
$$
k(\nu)\quad=\quad\prod_{B\in \nu}\left(\sum_{\delta\in \Pi(B)}\mm_B(\delta)\mu(\delta)\right)=\prod_{B\in \nu} k_B,
$$
which finishes the proof.
\hfill$\Box$\end{proof}

The inverse formula for moments in terms of cumulants follows directly by Proposition~\ref{prop:mobdual} and Lemma~\ref{lem:kisprod}. For every $I\subseteq [n]$ we have
\begin{equation}\label{eq:cumulinverse}
\mu_I\quad=\quad\sum_{\pi\in \Pi(I)}k(\pi)\quad=\quad\sum_{\pi\in \Pi(I)}\prod_{B\in \pi}k_B.
\end{equation}

We can directly generalize the definition of binary cumulants to \textit{binary $L$-cumulants}. Let  $L$ be a partition lattice of $[n]$. Define binary $L$-cumulants by 
\begin{equation}\label{eq:binaryell}
\ell_I\quad=\quad\sum_{\pi\in L(I)} \mm(\pi)\prod_{B\in \pi}\mu_B\qquad\mbox{for every }I\subseteq [n].
\end{equation}
By definition for every $I\subseteq [n]$ the maximal and minimal element of the lattice $L(I)$ coincide with the minimal and maximal element of $\Pi(I)$. In particular for every $L$ we have $\ell_i=\mu_i$ for $i=1,\ldots,n$; and $\ell_{ij}=\mu_{ij}-\mu_i\mu_j$ for all $1\leq i<j\leq n$. However, already when $n=3$ not all $L$-cumulants coincide with cumulants.
\begin{exmp}\label{ex:interval3}Let $n=3$ and consider $L$-cumulants induced by the lattice of interval partitions. The lattice $\cI([3])$ has four elements: $123$, $1|23$, $12|3$ and $1|2|3$ and $\mm(\pi)=(-1)^{|\pi|-1}$. Therefore, we have
$$
\ell_{123}\,\,=\,\,\mu_{123}-\mu_1\mu_{23}-\mu_{12}\mu_3+\mu_1\mu_2\mu_3.
$$
Compare this with the formula for $k_{123}$ in (\ref{eq:kappa3}) to note that not only the term $\mu_2\mu_{13}$ is missing now in the formula for $\ell_{123}$ but also the coefficient of $\mu_1\mu_2\mu_3$ is $1$ not $2$.
\end{exmp}

Let $\pi\in \Pi([n])$ be a set partition into blocks $B_1,\ldots,B_r$. Denote
$$\indep_{B\in \pi}X_B\quad:=\quad X_{B_1}\indep \cdots\indep X_{B_r}.$$ 
By (\ref{eq:indepinmoms}), $\indep_{B\in \pi}X_B$ if and only if  
\begin{equation}\label{eq:binaryindmoms}
\mu_I\,\,\,=\,\,\,\mu(\pi(I))\qquad\mbox{for every }I\subseteq [n],
\end{equation}
where $\pi(I)$ denotes $\pi$ constrained to elements in $B$. So for example the full independence is given by the minimal partition $\pi=1|2|\cdots|n$ and $\mu_I=\prod_{i\in I}\mu_i$ for every $I\subseteq [n]$. 

Below we list the basic facts about binary $L$-cumulants. They are proved in a more general setting in Section \ref{sec:XCproperties}. The following result implies that (P1) holds for binary $L$-cumulants.
\begin{prop}\label{prop:binindep}There exists a partition $\pi_0\in L$ such that $\indep_{B\in\pi_0}X_B$ if and only if $\ell(\pi)=0$ for all $\pi\not\leq\pi_0$, or equivalently, if $\ell_I=0$ unless $I$ is contained in one of the blocks of $\pi_0$ (equivalence follows from Theorem~\ref{th:P1}).
\end{prop}
\begin{proof}The 'if' part of the proposition is given in a more general setting in Proposition~\ref{prop:P1}.  To prove the opposite implication use Theorem~\ref{th:P1} to conclude that $\ell(\pi)=0$ for all $\pi\not\leq\pi_0$ implies that $\mu_I=\mu(\pi_0(I))$ for all $I\subseteq [n]$ which by (\ref{eq:binaryindmoms}) implies $\indep_{B\in \pi_0}X_B$.
\hfill$\Box$\end{proof}

\begin{exmp}\label{ex:n3interv} Consider the situation of Example~\ref{ex:interval3}, where $n=3$ and $L$-cumulants are defined by the lattice of interval partitions. If $X_1\indep (X_2,X_3)$ then $\mu_{123}=\mu_1\mu_{23}$, $\mu_{12}=\mu_1\mu_{2}$ and $\mu_{13}=\mu_1\mu_{3}$. It follows that $\ell_{12}=\ell_{13}=\ell_{123}=0$. On the other hand, the condition $X_2\indep (X_1,X_3)$ does not imply that $\ell_{123}=0$ because in this case
$$
\ell_{123}\,\,\,=\,\,\,\mu_2\mu_{13}-\mu_1\mu_2\mu_3,$$
which is zero only when in addition $\mu_{13}=\mu_1\mu_3$ and hence when $X_1\indep X_3$. Here there is no contradiction with Proposition~\ref{prop:binindep} because $2|13\notin \mathcal{I}([3])$. \end{exmp}

Under a minor additional condition the property (P2) also holds for binary $L$-cumulants.
\begin{prop}\label{prop:translbinary}Suppose that for every $i\in [n]$ the split $i|([n]\setminus i)$ lies in $L$. Let $\widetilde{X}=X+a$, where $a\in\R^{n}$ and, for every $I\subseteq[n]$, by $\widetilde{\ell}_I$ denote the corresponding $L$-cumulant of the subvector $X_I$. Then $\widetilde{\ell}_i=\ell_i+a_{i}$ for all $i=1,\ldots, n$ and  $\widetilde{\ell}_I= \ell_I$ for any $I\subseteq[m]$ such that $|I|\geq 2$. 
\end{prop}
\begin{proof}This follows from Proposition~\ref{prop:translation}.
\hfill$\Box$\end{proof}

Define \textit{central binary $L$-cumulants} by replacing moments $\mu_B$ in (\ref{eq:binaryell}) by central moments $\mu'_B$. For every $I\subseteq [n]$ the corresponding central binary $L$-cumulant is denoted by $\ell'_I$. 
\begin{lem}\label{lem:centrallcum}
Under the assumptions of Proposition~\ref{prop:translbinary} we have $\ell'_I=\ell_I$ for every $I\subseteq[n]$ such that $|I|\geq 2$.
\end{lem}
\begin{proof}
Central binary $L$-cumulants of $X$ can be alternatively defined as binary $L$-cumulants of $\widetilde{X}$, where $\widetilde{X}_i=X_i-\E X_i$. The lemma follows from Proposition~\ref{prop:translbinary}.
\hfill$\Box$\end{proof}

In the next section we show how all these ideas can be applied in algebraic geometry.

\subsection{Geometric applications}\label{sec:geometry}

We consider algebraic varieties in either the real space $\R^{2^n}=\R^{2\times\cdots\times 2}$ or its complexification $\C^{2^n}=\C^{2\times\cdots\times 2}$, or projectivization $\P^{2^n-1}=\P(\C^{2\times\cdots\times 2})$. Each component $\C^2$ (or $\R^2$) has basis $e_0,e_1$ so that $e_{i_1}\otimes \cdots\otimes e_{i_n}$ corresponds to $I\subseteq [n]$ for $i_j=1$ if $j\in I$ and $i_j=0$ otherwise. For example, if $n=2$ and $\mu\in \C^{2\times 2}$ then we write $\mu$ in our basis as
$$
\mu\quad =\quad\mu_\emptyset \,e_0\otimes e_0+\mu_1\, e_1\otimes e_0+\mu_2 \,e_0\otimes e_1+\mu_{12}\, e_1\otimes e_1.
$$
Formula (\ref{eq:binaryell}) gives an isomorphism of the affine subspace $\mu_\emptyset=1$ in $\R^{2^n}$ (or $\C^{2^n}$), which forms a Zariski open subset of $\P^{2^n-1}$. The inverse map is computed in a more general case in (\ref{eq:inversexi}).

We first show that some basic operations on the random vector $X$ can encode interesting actions on the space of $2\times\cdots\times 2$ tensors. Define $\widetilde{X}$ such that $\widetilde{X}_i=\lambda_iX_i$ for $\lambda_i\in \C\setminus\{0\}$ for $i=1,\ldots,n$.  Multiplying each $X_i$ by $\lambda_i$ results in the change of moments from $\mu_I$ to $\widetilde{\mu}_I=\prod_{i\in I}\lambda_i\mu_{I}$ and hence it corresponds to the action of the group $D^n$, where $D$ a group of diagonal matrices of the form 
$$
\left[\begin{array}{cc}
1 & 0\\
0 & \lambda
\end{array}
\right]\qquad \mbox{for }\lambda\in \mathbb{C}\setminus\{0\}.
$$ 
Because $L$-cumulants are multilinear functions of the moments we conclude that this action is the same on the level of $L$-cumulants. We have $\widetilde{\ell}_I=\prod_{i\in I}\lambda_i\ell_{I}$ for every $I\subseteq[n]$.

Suppose now that $\widetilde{X}=X+b$, for $b=(b_1,\ldots,b_n)\in \C^n$, and consider the group $U(2)^n$ where $U(2)$ is the unipotent group
of $2\times 2$-matrices of the form
$$
\left[\begin{array}{cc}
1 & 0\\
\lambda & 1
\end{array}
\right]\qquad \mbox{for }\lambda\in \mathbb{C}.
$$ 
Adding $b$ to the vector $X$ corresponds to the action of $U(2)^n$, with $\lambda_i=b_i$ for $i=1,\ldots,n$, on the space of moments. We illustrate this with an example that easily generalizes.
\begin{exmp}\label{ex:n2tensrs}
Let $n=2$ and denote by $\widetilde{\mu}=[\widetilde{\mu}_I]$ the moments of the vector $\widetilde{X}=X+b$. We have $\widetilde{\mu}_\emptyset=1$, $\widetilde{\mu}_i=\E(X_i+b_i)=\mu_i+b_i$ for $i=1,2$ and
$$
\widetilde{\mu}_{12}\,\,:=\,\,\E[(X_1+b_1)(X_2+b_2)]\,\,=\,\,\mu_{12}+b_1\mu_2+\mu_1b_2+b_1b_2.
$$
Write $\mu=[\mu_I]\in \C^{2\times 2}$: 
$$\mu\quad=\quad e_0\otimes e_0+\mu_1 e_1\otimes e_0+\mu_2 e_0\otimes e_1+\mu_{12}e_1\otimes e_1.$$
After applying the action of $U(2)^2$ with $\lambda_i=b_i$ for $i=1,2$ we obtain
\begin{eqnarray*}
\widetilde{\mu}\,\,\,&=&\,\,\, (e_0+b_1e_1)\otimes (e_0+b_2e_1)+\mu_1 e_1\otimes (e_0+b_2e_1)+\mu_2 (e_0+b_1e_1)\otimes e_1+\mu_{12}e_1\otimes e_1\\
&=\,\,\,& e_0\otimes e_0+(\mu_1+b_1) e_1\otimes e_0+(\mu_2+b_2) e_0\otimes e_1+(\mu_{12}+b_1\mu_2+\mu_1b_2+b_1b_2)e_1\otimes e_1\\
&=\,\,\,& e_0\otimes e_0+\widetilde{\mu}_1 e_1\otimes e_0+\widetilde{\mu}_2 e_0\otimes e_1+\widetilde{\mu}_{12}e_1\otimes e_1,
\end{eqnarray*}
which confirms that translating $X$ by $b\in \R^n$ corresponds to the action of $U(2)^n$ on $\mu$.
\end{exmp}

For every $I\subseteq[n]$, denote by $\widetilde{\ell}_I$ the $L$-cumulant of $\widetilde{X}_I$. By Proposition~\ref{prop:translbinary}, whenever every split $i|([n]\setminus i)$ lies in $L$, this complicated transformation of moments induced by $U(2)^n$ translates to a very simple transformation of cumulants. We have $\widetilde{\ell}_i=\ell_i+b_i$ for $i\in [n]$ and   $\widetilde{\ell}_I=\ell_I$ for all $I\subseteq[n]$ such that $|I|\geq 2$ and hence all the higher order $L$-cumulants are invariant with respect to the action of $U(2)^n$ on the space of moments. 

Changing values of the binary
variables $X_i$ from $0,1$ to $b_i, a_i$, means defining a new random vector $\widetilde{X}$ such that $\widetilde{X}_i=(a_i-b_i)X_i+b_i$. We have just shown that changing values of the components of $X$ corresponds to  a natural action of the $n$-dimensional torus $(\mathbb{C}^*)^n$ with coordinates $a_i-b_i$ on the space  $\mathbb{C}^{2^n-n-1}$ whose coordinates are the higher order $L$-cumulants $\ell_I$, $|I| \geq 2$. More specifically the $L$-cumulants of $\widetilde{X}$, such that $\widetilde{X}_i=(a_i-b_i)X_i+b_i$, are transformed by 
\begin{equation*}
\widetilde{\ell}_{I}\,\,\, = \,\,\, \ell_I \cdot \prod_{i\in I}(a_i-b_i) \qquad\mbox{ for all } I\subseteq [n] 
\hbox{ and } |I|\geq 2
\end{equation*}
and $\widetilde{\ell}_{i}=(a_{i}-b_{i})\ell_{i}+b_{i}$ for $i=1,\ldots, n$. This leads to the following result.
\begin{thm} \label{thm:invariant} Suppose that for every $i\in [n]$ the split $i|([n]\setminus i)$ lies in $L$. Then a subvariety of $\mathbb{C}^{2^n-1}$ is
invariant under changing values of components of $X$
if and only it is defined by  $\mathbb{Z}^n$-homogeneous polynomials
in $\ell_I$ with~$|I| \geq 2$.
\end{thm}
\begin{proof}
See the proof of \cite[Theorem~3.1]{pwz-2011-bincum}.
\hfill$\Box$\end{proof}

Note that if a variety is invariant under the action of the \textit{special linear group ${\rm SL}(2)^n$} then in particular it is invariant under $U(2)^n$. 
\begin{cor} \label{cor:invariant}
Suppose that $L$ is a partition lattice of $[n]$ such that for every $i\in [n]$ the split $i|([n]\setminus i)$ lies in $L$. Let $V$ be a subvariety of the affine open subset given by $\mu_\emptyset = 1$ in 
the projective space $\mathbb{P}(\mathbb{C}^{2 \times  \cdots \times 2})$
and let $\overline{V}$ denote its closure in that projective space.
If  $\overline{V}$ is invariant
under the action of $SL(2)^n$ then the ideal $I_V$ that defines $V$ is generated by
$\mathbb{Z}^n$-homogeneous polynomials in the $L$-cumulants $\ell_I$ with $|I| \geq 2$.
\end{cor}

Another important reason why $L$-cumulants may be useful, apart from their invariance properties, is related to property (P1). Denote by ${\rm Seg}((\P^1)^n)$ the \textit{Segre variety}, which is an embedding of $(\P^1)^n$ into $\P^{2^n-1}$. In statistics the Segre variety corresponds to the full independence model $X_1\indep \cdots\indep X_n$. In particular Proposition \ref{prop:binindep} implies that the image of ${\rm Seg}((\P^1)^n)$ in the space given by $L$-cumulants is an affine subspace given by $\ell_I=0$ for all $|I|\geq 2$ (see also \cite[Remark 3.4]{pwz-2011-bincum}). Moreover, $L$-cumulants seem to be helpful also in the analysis of other algebraic varieties related to the Segre variety ${\rm Seg}((\P^1)^n)$. For example the tangential variety ${\rm Tan}((\P^1)^n)$ is toric when expressed in cumulants (see \cite[Theorem~4.1]{pwz-2011-bincum}). In the following section we show how $L$-cumulants defined by a tree partition lattice can help to analyze the \textit{secant variety ${\rm Sec}((\P^1)^n)$}.

\subsection{Binary tree cumulants for secant varieties}\label{sec:tressforsecs}

In \cite{pwz-2010-identifiability} we defined tree cumulants, which gave a better understanding of certain statistical models related to trees. We write more on that in Section \ref{sec:treecum}. In this section we show how tree cumulants can be used to study secant varieties. Recall from Definition~\ref{def:partitions} that, for a fixed tree $T$ with $n$-leaves, $\cT^{T}([n])$ denotes the lattice of tree partitions of $[n]$ induced by $T$. Moreover, $\cT^{T}(I)$ is the lattice of all tree partitions of $I$ induced by $T(I)$, which is the smallest subtree of $T$ containing all leaves in $I$. The tree cumulant of the subvector $X_I$ for every $I\subseteq [n]$ is denoted by $\mathfrak{t}_I$. Tree cumulants are $L$-cumulants and hence defined by (\ref{eq:binaryell}):
\begin{equation}\label{eq:treecumuls}
\mathfrak{t}_I\quad =\quad \sum_{\pi\in \cT^{T}(I)}\mm(\pi)\prod_{B\in \pi}\mu_B,\qquad\mbox{for all }I\subseteq [n].
\end{equation}

\begin{rem}\label{rem:centraltree}In \cite[Section 3.2]{pwz-2010-identifiability} binary tree cumulants were defined in terms of central moments by  
$$
\widetilde{\mathfrak{t}}_I\,\,\,=\,\,\,\sum_{\pi\in \cT^{T}(I)}\mm(\pi)\prod_{B\in \pi}\mu'_B\qquad\mbox{for all }I\subseteq [n], \,|I|\geq 2,
$$
and $\widetilde{\mathfrak{t}}_i=\mu_i$ for $i\in [n]$. In particular  $\widetilde{\mathfrak{t}}_I$ for all $|I|\geq 2$ is just the corresponding central $L$-cumulant. Let $i\in[n]$ be one of the leaves. Removing the edge incident with $i$ induces a split $i|([n]\setminus i)$ and hence the assumption of Proposition~\ref{prop:translbinary} holds and, by Lemma~\ref{lem:centrallcum}, it follows that $\widetilde{\mathfrak{t}}_I=\mathfrak{t}_I$ for all $I\subseteq [n]$. In particular, both the definition in \cite{pwz-2010-identifiability} and the one given in (\ref{eq:treecumuls}) are equivalent.\end{rem}

Let $L$ be the lattice of tree partitions induced by the caterpillar tree in Figure \ref{fig:cater}. For example if $n=4$ then the induced lattice is given in Figure \ref{fig:treeposet4}. We first show how to compute $L$-cumulants $[\mathfrak{t}_I]$ without  computing the M\"{o}bius function on the lattice $L$. By Remark \ref{rem:centraltree} we can replace moments by central moments in the formula for $\mathfrak{t}_I$ for all $I\subseteq[n]$ such that $|I|\geq 2$. This is very convenient because $\prod_{B\in \pi}\mu'_B$ is zero whenever $\pi$ contains a singleton block. Note that the elements of $L$ with no singleton blocks correspond to all interval partitions with no singleton blocks. If $n=4$ then the elements of $L$ with no singleton blocks are the two boldfaced elements in Figure \ref{fig:treeposet4}. This gives that for all $I\subseteq[n]$ such that $|I|\geq 2$:
$$
\mathfrak{t}_{I}\quad=\quad\sum_{\pi\in L(I)}\mm({\pi})\prod_{B\in \pi} \mu'_B\quad=\quad\sum_{\pi\in \cI(I)}\mm({\pi})\prod_{B\in \pi} \mu'_B.
$$
Both sums above are over all partitions in a poset of all interval partitions with no singleton blocks. Hence, both M\"{o}bius functions constrained to this poset need to coincide. The gain is that we already computed the M\"{o}bius function on the right-hand side explicitly obtaining $\mm(\pi)=(-1)^{|\pi|-1}$ (see the end of Section \ref{sec:XClattice}).

This allows us to write the map from moments $[\mu_I]$ to  tree cumulants $[\mathfrak{t}_I]$ of the caterpillar tree as a composition of two maps: from moments to central moments and from central moments to tree cumulants induced by the caterpillar tree. We will show in the end of Section~\ref{sec:props} that the first map can be written as
\begin{equation*}
\mu_I'\quad =\quad \sum_{B\subseteq I}(-1)^{|I\setminus B|}\mu_B\prod_{i\in I\setminus B}\mu_i \qquad\mbox{for all }I\subseteq [n],\, |I|\geq 2,
\end{equation*}
and we have just shown that the second map is given by $\mathfrak{t}_i=\mu_i$ for $i=1,\ldots, n$, and 
\begin{equation*}
\mathfrak{t}_{I}\quad=\quad\sum_{\pi\in \cI(I)}(-1)^{|\pi|-1}\prod_{B\in \pi} \mu'_B\qquad\mbox{for all }|I|\geq 2.
\end{equation*}
In particular, if $n=4$ then $\mathfrak{t}_I=\mu'_I$ for all $2\leq |I|\leq 3$ and
$$
\mathfrak{t}_{1234}\,\,\,=\,\,\,\mu'_{1234}-\mu'_{12}\mu_{34}'.
$$

\begin{figure}[htp!]
\centering
\tikzstyle{vertex}=[circle,fill=black,minimum size=5pt,inner sep=0pt]
\tikzstyle{hidden}=[circle,draw,minimum size=5pt,inner sep=0pt]
  \begin{tikzpicture}
  \node[vertex] (1) at (0,1)  [label=above:$1$] {};
    \node[vertex] (2) at (1,1)  [label=above:$2$] {};
    \node[vertex] (3) at (2,1) [label=above:$3$]{};
    \node[vertex] (4) at (4,1) [label=above:$n$]{};
    \node[hidden] (a) at (0,0)  {};
    \node[hidden] (b) at (1,0) {};
    \node[hidden] (c) at (2,0) {};
     \node (d) at (3,0)  {$\cdots$};
     \node[hidden] (e) at (4,0)  {};
           \draw (a) to (b);
           \draw (c) to (b);
           \draw (c) to (d);
           \draw (e) to (d);
               \draw (a) to (1);
    \draw (b) to (2);
    \draw (c) to (3);
    \draw (e) to (4);
  \end{tikzpicture}\caption{A caterpillar tree with $n$ leaves/legs.}\label{fig:cater}
\end{figure}
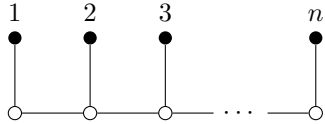
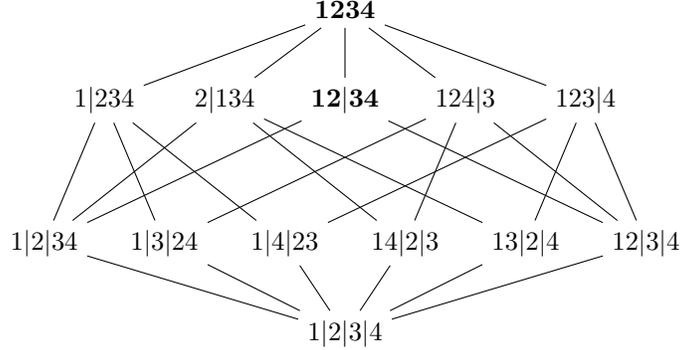
\begin{figure}[htp!]
\centering
\tikzstyle{vertex}=[circle,fill=black,minimum size=5pt,inner sep=0pt]
\tikzstyle{hidden}=[circle,draw,minimum size=5pt,inner sep=0pt]
  \begin{tikzpicture}
    \node (1234) at (0,5)   {$\mathbf{1234}$};
    \node (a1) at (-3.2,3.8)   {$1|234$};
    \node (a2) at (-1.6,3.8) {$2|134$};
    \node (a3) at (0,3.8) {$\mathbf{12|34}$};
    \node (a4) at (1.6,3.8)  {$124|3$};
    \node (a5) at (3.2,3.8) {$123|4$};
    \node (b1) at (-4,1.9) {$1|2|34$};
    \node (b2) at (-2.4,1.9)  {$1|3|24$};
    \node (b3) at (-0.8,1.9)  {$1|4|23$};
    \node (b4) at (.8,1.9)  {$14|2|3$};
    \node (b5) at (2.4,1.9)  {$13|2|4$};
    \node (b6) at (4,1.9)  {$12|3|4$};
    \node (c) at (0,0.7)  {$1|2|3|4$};
 \draw (1234) to (a1);
 \draw (1234) to (a2);
 \draw (1234) to (a3);
 \draw (1234) to (a4);
 \draw (1234) to (a5);
 \draw (a1) to (b1);
 \draw (a1) to (b2);
 \draw (a1) to (b3);
 \draw (a2) to (b1);
 \draw (a2) to (b4);
 \draw (a2) to (b5);
 \draw (a3) to (b1);
 \draw (a3) to (b6);
 \draw (a4) to (b2);
 \draw (a4) to (b4);
 \draw (a4) to (b6);
 \draw (a5) to (b3);
 \draw (a5) to (b5);
 \draw (a5) to (b6);
 \draw (b1) to (c);
 \draw (b2) to (c);
 \draw (b3) to (c);
 \draw (b4) to (c);
 \draw (b5) to (c);
 \draw (b6) to (c);
  \end{tikzpicture}\caption{The Hasse diagram of the lattice of tree partitions induced by the tree in Figure \ref{fig:cater} if $n=4$.}\label{fig:treeposet4}
\end{figure}

We use this new coordinate system to study the secant variety  ${\rm Sec}((\P^1)^n)$. As an example consider the case when $n=4$.
\begin{exmp}\label{ex:secant}
The secant variety ${\rm Sec}((\P^1)^4)$ is a projective variety in $\P^{15}$ parametrized by $9$ copies of $\P^1$ with coordinates $({t}_0,t)$,  $({a}_{0i},a_i)$ and $({b}_{0i},b_i)$ for $i=1,2,3,4$. The parametrization is given by
$$
\mu_I\quad=\quad t_0\prod_{i\in I^c} {a}_{0i}\prod_{i\in I} a_i\,\,\,+\,\,\,{t}\prod_{i\in I^c} {b}_{0i}\prod_{i\in I} b_i\qquad\mbox{for all }I\subseteq[4],
$$
where $I^c$ denotes the complement of $I$ in $\{1,2,3,4\}$ and $\mu=[\mu_I]$ denotes the coordinates of the projective space $\P^{15}$. We want to describe the image of an open subset of the parameter space given by ${a}_{0i}={b}_{0i}=1$ for $i\in \{1,2,3,4\}$ and $t_0=1-t$. This image is described by
\begin{equation}\label{eq:paramnclaw}
\mu_I\quad=\quad (1-t)\prod_{i\in I} a_i\,\,\,+\,\,\,t\prod_{i\in I} b_i
\end{equation}
and in particular $\mu_\emptyset=1$. 

Earlier in this section we explained how to compute $[\mathfrak{t}_I]$ from moments as a composition of two simple maps. From this we can also compute the induced parametrization directly. Here we will show an alternative way of proceeding for the secant variety ${\rm Sec}((\P^1)^4)$ to present some other available techniques. First, use the parametrization of the secant in terms of classical cumulants. This parametrization was given in  \cite[Equations (18) and (19)]{pwz-2011-bincum}, which implies that for every $i<j<k$
 \begin{equation}
 \label{eq:kappa}
 \begin{matrix} 
  k_{ij} &=& t(1-t)(b_i-a_i)(b_j-a_j)  \\
  k_{ijk} &=& t(1-t)(1-2t)(b_i-a_i)(b_j-a_j)(b_k-a_k)  \\
  k_{1234} &=& t(1-t)(6t^2-6t+1)\prod_{i=1}^4 (b_i-a_i)  .
 \end{matrix}
\end{equation}
Now we change coordinates from cumulants to binary tree cumulants $[\mathfrak{t}_I]$ using Proposition~\ref{prop:LinK}. In particular, as explained in Example~\ref{ex:notinPi}, since $13|24$ and $14|23$ are the only partitions in $\Pi([4])$ which are not tree partitions of the caterpillar tree in Figure \ref{fig:cater} for $n=4$, this yields
\begin{equation}\label{eq:l1234bis}
\mathfrak{t}_{1234}\,\,\,=\,\,\,k_{1234}+k_{13}k_{24}+k_{14}k_{23}
\end{equation}
and $\mathfrak{t}_I=k_I$ for all $I\subseteq [4]$ such that $|I|\leq 3$. From this it follows that for every $I\subseteq \{1,2,3,4\}$ such that $|I|\geq 2$:
\begin{equation}\label{eq:secantparint}
\mathfrak{t}_I\,\,\,=\,\,\,t(1-t)(1-2t)^{|I|-2}\prod_{i\in I} (b_i-a_i),
\end{equation}
which for $\mathfrak{t}_{1234}$ can be verified by direct computations. Now we can immediately check that 
$$\mathfrak{t}_{I\cup J} \mathfrak{t}_{I'\cup J'}\,\,\,-\,\,\, \mathfrak{t}_{I\cup J'} \mathfrak{t}_{I'\cup J}\quad =\quad 0
$$
holds on ${\rm Sec}((\P^1)^4)$ for all distinct $I,I'\in \{\{i\},\{j\},\{i,j\}\}$ and $J,J'\in \{\{k\},\{l\},\{k,l\}\}$ and every split $ij|kl$ of $\{1,2,3,4\}$. For example $12|34$ leads to a set of equations including $\mathfrak{t}_{13} \mathfrak{t}_{24}- \mathfrak{t}_{14} \mathfrak{t}_{23}=0$ and $\mathfrak{t}_{1234} \mathfrak{t}_{13}- \mathfrak{t}_{123} \mathfrak{t}_{134}=0$. \end{exmp}

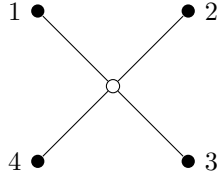
\begin{figure}[htp!]
\centering
\tikzstyle{vertex}=[circle,fill=black,minimum size=5pt,inner sep=0pt]
\tikzstyle{hidden}=[circle,draw,minimum size=5pt,inner sep=0pt]
  \begin{tikzpicture}
  \node[vertex] (1) at (0,1)  [label=left:$1$] {};
    \node[vertex] (2) at (2,1)  [label=right:$2$] {};
    \node[vertex] (3) at (2,-1) [label=right:$3$]{};
    \node[vertex] (4) at (0,-1) [label=left:$4$]{};
    \node[hidden] (a) at (1,0)  {};
               \draw (a) to (1);
    \draw (a) to (2);
    \draw (a) to (3);
    \draw (a) to (4);
  \end{tikzpicture}\caption{A 4-star tree.}\label{fig:4star}
\end{figure}

This simple example can be generalized using the link between the secant varieties and certain statistical models (see \cite[Section 4.1]{oberwolfach2009}). Define for any two disjoint $A,C\subseteq [n]$ the conditional probability of $X_A$ given $X_C$ as:
\begin{equation*}
p_{A|C}(x_A|x_C)\,\,\,:=\,\,\,\frac{p_{A\cup C}(x_A,x_C)}{p_C(x_C)} \quad\mbox{for all }x_C\in \cX_C \mbox{ s.t. } p_C(x_C)\neq 0.
\end{equation*}
For any function $f$ of $X_A$ define the \textit{conditional expectation of $f(X_A)$ given $X_C$} as a function of $X_C$ given for any $x_C\in \cX_C$ 
$$
\E[f(X_A)|X_C=x_C]\,\,\,=\,\,\,\sum_{x_A\in \cX_A}p_{A|C}(x_A|x_C)f(x_A).
$$
We denote this conditional expectation by $\E[f(X_A)|X_C]$. If $f(X_A)=\prod_{i\in A}X_i$ then we simply write $\mu_A^C$ and $\mu_A^C(x_C)=\E[\prod_{i\in A}X_i|X_C=x_C]$. Note that $\mu_A^C$ is a random variable itself.  

Similarly as in the case of Lemma \ref{lem:indinmom} we can show that for disjoint $C,B_1,\ldots,B_r\subseteq [n]$ the $X_{B_i}$'s are \textit{jointly independent given $X_C$} if 
\begin{equation*}
\mu_{A_1\cdots A_r}^C\,\,\,=\,\,\,\prod_{i=1}^r \mu_{A_i}^C\qquad\mbox{for all }A_i\subseteq B_i,\, i=1,\ldots,r,
\end{equation*}
In this case the marginal distribution of $X_C$ satisfies
\begin{equation}\label{eq:likesecant}
\mu_{A_1\cdots A_r}\,\,\,=\,\,\,\E[\mu_{A_1\cdots A_r}^C]\,\,\,=\,\,\,\sum_{x_C\in \cX_C}p_C(x_C)\mu_{A_1\cdots A_r}^C(x_C).
\end{equation}

For a statistician the parametrization in (\ref{eq:paramnclaw}) corresponds to the parametrization of moments of the binary $4$-star tree model (naive Bayes model) as given in Figure \ref{fig:4star}. The leaves of this tree correspond to a vector $X=(X_1,X_2,X_3,X_4)$ of binary \textit{observed} variables and the inner node corresponds to a binary variable $Y$ which is not observed. This model contains all possible moments of a binary vector $X$ such that all components of $X$ are jointly independent given $Y$. The parametrization in  (\ref{eq:paramnclaw}) is a special version of (\ref{eq:likesecant}).

The fact that (\ref{eq:paramnclaw}) can be rewritten in the easier form in (\ref{eq:secantparint}) for any $n\geq 4$ follows from more general considerations in \cite[Section 4]{pwz-2010-identifiability}. We obtain the following procedure:
\begin{itemize}
\item[1.] Consider any trivalent tree with $n$ leaves, that is a tree such that each inner node has valency three.
\item[2.] Compute tree cumulants induced by this trivalent tree.
\item[3.] The induced parametrization of the $n$-star tree model in the coordinate system constructed in step 2 is  (\ref{eq:secantparint}), where now $I\subseteq [n]$ for $n\geq 4$. For more details check Section \ref{sec:treecum}.
\end{itemize}
Of course, since we can pick any trivalent tree in step 1, the most natural choice is to pick the caterpillar tree. This is mainly because the computation of the corresponding tree cumulants is simple as it was presented earlier in this section. Now from the parametrization in (\ref{eq:secantparint}) we easily verify that 
\begin{equation}\label{eq:eqssecant}
\mathfrak{t}_{I\cup J} \mathfrak{t}_{I'\cup J'}- \mathfrak{t}_{I\cup J'} \mathfrak{t}_{I'\cup J}\quad =\quad 0
\end{equation}
holds on ${\rm Sec}((\P^1)^n)$ for all non-empty subsets $I,I'\subseteq A$ and $J,J'\subseteq B$ where $A|B$ is a split of $[n]$.

\begin{rem}\label{rem:naively}
It may seem that a more natural way to proceed in Example~\ref{ex:secant} was to construct tree cumulants induced directly by partitions of the $4$-star tree in Figure \ref{fig:4star}. The tree partitions of the $4$-star tree are equal to one-cluster partitions from Definition~\ref{def:partitions}. By Proposition \ref{prop:centralm} this partition lattice induces central moments $\mu'_I$. To compute the induced parametrization of the central moments note that $\mu'_I=k_I$ for all  $2\leq |I|\leq 3$. A direct check shows that 
$$
\mu'_{1234}\,\,=\,\,t(1-t)(3t^2-3t+1)\prod_{i=1}^4 (b_i-a_i)
$$
and we find that the relation between $\mu'_{1234}$ and other central moments is more complicated than in the case of  tree cumulants induced by the caterpillar tree.  In particular, the corresponding equations are no longer binomial like in (\ref{eq:eqssecant}). 
\end{rem}

\section{The definition of $L$-cumulants}\label{sec:XCdef}

Let $A=\{i_1,\ldots,i_d\}$ be a multiset. We define its \textit{multisubset} $B\subseteq A$ as a multiset $B=\{i_j:\,j\in I\}$ for some $I\subseteq [d]$. For example if $A=\{1,1,2,2\}$ then $A$ has, among others, four multisubsets of the form $\{1,2\}$. Let $X$ be a finite discrete random vector with values in $\cX$ and let $\cA(\cX)$ be the family of multisets associated to $\cX$ as given in (\ref{eq:defAX}). Consider any family $\mathbf{L}=(L({A}))_{A\in \cA(\cX)}$ of partition lattices such that $L(A)$ is a subposet of $\Pi(A)$ for every $A\in \cA(\cX)$. Assume that the maximal and minimal elements of $L(A)$ coincide with the maximal and the minimal element of $\Pi(A)$ and denote them by ${A}$ and $\hat{0}_{A}$ respectively. Moreover, for every $B\subseteq A$ the map $L(A)\rightarrow L(B)$ are surjections given by constraining partitions of $A$ to $B$.  Note that in particular,  $L(A)$ need not be a sublattice of $\Pi(A)$ because the join and the meet operator of $L(A)$ and $\Pi(A)$ may differ.

 The first two trivial examples of a family $\mathbf{L}$ as above is $\mathbf{\Pi}=(\Pi(A))_{A\in \cA(\cX)}$ and $\mathbf{L}$ such that for every $A\in \cA(\cX)$, $|A|\geq 2$, the lattice $L(A)$ is given by just two elements $\hat{0}_A$ and $A$. Other interesting examples are obtained from Definition~\ref{def:partitions} (excluding tree partitions), where $L(A)$ is assumed to be isomorphic to $L(|A|)$. The corresponding families of lattices are denoted by $\mathbf{NC}$ (non-crossing), $\mathbf{I}$ (interval) and $\mathbf{C}$ (one-cluster). A definition of tree cumulants in this case requires construction of an $\hat{A}$-labelled tree, where $\hat{A}$ is  the maximal multiset $A$ in $\cA(\cX)$ corresponding to $x=(r_1-1,\ldots,r_n-1)$. This construction is not unique and for that reason we discuss tree cumulants only in very concrete examples.
 
By $\mm_{A}$ we denote the M\"{o}bius function on $ L(A)$. The lattice will be always obvious from the context so we omit it in the notation. When $A$ is also clear from the context we just write $\mm$. 
\begin{defn}[$L$-cumulants]\label{def:xicumulants}
Let $X=(X_{1},\ldots, X_{n})$ be a random vector. For any $A\in \cA(\cX)$  and $\nu\in L(A)$ define 
\begin{equation}\label{xi-cumul}
\ell(\nu)\quad=\quad\sum_{\pi\leq \nu}\mm_{A}(\pi,\nu)\mu(\pi),
\end{equation}
where $\mu(\pi)=\prod_{B\in \pi} \mu_{B}$. Then $\ell_{A}:=\ell({A})$ is the \textit{$L$-cumulant of $X_{A}$}. 
\end{defn}

If $\mathbf{L}=\mathbf{\Pi}$ then, because $\Pi(A)\simeq \Pi({[|A|]})$, we obtain the formula in (\ref{eq:cumA}) and hence this definition generalizes the classical cumulants. Other known $L$-cumulants were defined in the non-commutative probability literature. These are $L$-cumulants defined by $\mathbf{NC}$ and $\mathbf{I}$, which  are called \textit{free cumulants} and \textit{Boolean cumulants} respectively (see \cite{speicher_freeprob97,speicher97}).

The map (\ref{xi-cumul}) is invertible with the inverse given by the M\"{o}bius inversion formula in Proposition~\ref{prop:mobdual}. Thus for every $A\in \cA(\cX)$
\begin{equation}\label{eq:inversexi}
\mu_{A}\,\,\,=\,\,\,\mu({A})\,\,\,=\,\,\,\sum_{\pi\in  L(A)}\ell(\pi).
\end{equation}
Note that in general $\ell(\pi)\neq \prod_{B\in \pi}\ell_B$, as it was the case for cumulants. However, $\ell(\pi)= \prod_{B\in \pi}\ell_B$ whenever $\mathbf{L}$ satisfies the following condition:
\begin{itemize}
\item[(C0)] For every $A\in\cA(\cX)$ and for any two partitions $\pi, \nu\in L(A)$ the interval $[\pi,\nu]$ is isomorphic to  the product of intervals $\prod_{B\in \nu} [\pi(B),{B}]\subseteq \prod_{B\in \nu} L(B)$.
\end{itemize}
Condition (C0) is not very restrictive. In fact all partition lattices mentioned in Definition~\ref{def:partitions} satisfy this property. If (C0) holds, then, by Proposition~\ref{prop:mobius-product}, the M\"{o}bius function on  $L(A)$ satisfies $\mm_{A}(\pi,\nu)=\prod_{B\in \nu}\mm_{B}(\pi({B}),{B})$. In particular (\ref{eq:inversexi}) becomes
\begin{equation*}
\mu_{A}\,\,\,=\,\,\,\sum_{\pi\in  L(A)}\prod_{B\in \pi}\ell_B,
\end{equation*}
and the proof of this follows essentially the proof of Lemma~\ref{lem:kisprod}. 

\begin{rem}By the moment aliasing there is a one-to-one correspondence between the probabilities $P=[p(x)]_{x\in \cX}$ and moments $M=[\mu_A]_{A\in \cA(\cX)}$ and hence also $L$-cumulants $\mathcal{L}=[\ell_A]_{A\in \cA(\cX)}$.\end{rem}

Unlike in the case of cumulants, for general $L$-cumulants no generating function is known. It may be then useful to realize that $L$-cumulants can be expressed in terms of classical cumulants in a rather simple manner. The following result generalizes Theorem~4.1 in \cite{lehner2002free}.
\begin{prop}\label{prop:LinK}
Let $L(A)$  be a lattice of set partitions of $A$ in the family $\mathbf{L}$ and let $\Pi^{*}$ denote the set of elements $\pi\in \Pi(A)$ such that $[\pi,{A}]\cap L(A)=\{A\}$, where the interval $[\pi,A]$ is taken in $\Pi(A)$. We have
$$
\ell_{A}\,\,\,=\,\,\,\sum_{\pi\in\Pi^{*}}k(\pi)\,\,\,=\,\,\,\sum_{\pi\in\Pi^{*}}\prod_{B\in \pi}k_B.
$$
\end{prop}
\begin{proof}
In this proof $\delta\leq_{\Pi}\pi$ means that $\delta\leq \pi$ and $\delta\in \Pi(A)$. Similarly $\pi\geq_{L}\delta$ denotes $\pi\geq \delta$ and $\pi\in L(A)$. Expressing the $L$-cumulant in terms of moments and then the moments in terms of classical cumulants gives
\begin{eqnarray*}
\ell_{A}&=&\sum_{\pi\in L(A)}\mm(\pi)\prod_{B\in \pi}\left(\sum_{\delta_B\in \Pi(B)} \prod_{C\in \delta_B}k_C\right)=\\
&=&\sum_{\pi\in L(A)}\mm(\pi)\sum_{\delta\leq_{\Pi} \pi} \prod_{B\in \delta}k_B.
\end{eqnarray*}
For every $\delta\in \Pi(A)$ let $\bar{\delta}$ denote the smallest element of $L(A)$ such that $\delta\leq_{\Pi} \bar{\delta}$. Then, by changing the order of summation, the above equation can be rewritten as
$$
\ell_{A}\,\,\,=\,\,\,\sum_{\delta\in \Pi(A)}\prod_{B\in \delta}k_B\left(\sum_{\pi\geq_{L}\bar{\delta}}\mm(\pi)\right).
$$
By (\ref{eq:sumintzero}) the sum in brackets vanishes whenever $\bar{\delta}\neq {A}$. Therefore the whole expression is equal to 
$\sum_{\delta\in \Pi^{*}}\prod_{B\in \delta}k_B$.
\hfill$\Box$\end{proof}

\begin{exmp}\label{ex:notinPi}Let $n=4$ and let $L$ be the lattice of all set partitions in Figure \ref{fig:treeposet4}. The only partitions of $\Pi([4])$ which are not in $L$ are $13|24$ and $14|23$. Hence, they are also the only partitions satisfying the condition $[\pi,[4]]\cap L=[4]$. This, by Proposition \ref{prop:LinK}, gives the formula for $\mathfrak{t}_{1234}$ given in (\ref{eq:l1234bis}).
\end{exmp}

\section{Basic properties of $ L$-cumulants}\label{sec:XCproperties}

In this section we show that $L$-cumulants satisfy properties similar to (P1)-(P4). The following lemma is central to most of the proofs of this section. It was first formulated by Weisner \cite{weisner35} in a special case and then generalized by Rota \cite{rota1964fct} for general lattices (see the corollary on page 351 therein). 
\begin{lem}\label{lem:sumzeroM}\label{lem:weisner}
Let $L$ be a finite lattice with at least two elements, and let $\pi_{0}\in L$ be such that $\pi_{0}\neq \hat{1}$. Then for any $\delta\in L$
$$
\sum_{\pi:\pi\wedge \pi_{0}=\delta}\mm(\pi)\,\,\,=\,\,\,0.
$$
\end{lem}
A special case of this result, when $\delta=\hat{0}$, is given in \cite[Corollary 3.9.3]{stanley2006enumerative}. It is a useful exercise to see that the proof given there generalizes to provide a proof of Lemma \ref{lem:weisner}.

\subsection{Independence and semi-invariance}\label{sec:props}

To show that property (P1) holds for $L$-cumulants we first prove a more algebraic version of this result. This result is directly linked to the definition of independence formulated in terms of moments in (\ref{eq:indepinmoms}).
\begin{thm}\label{th:P1}
Consider the $ L$-cumulant of $X=(X_{1},\ldots, X_{n})$ as in Definition~\ref{def:xicumulants}. The following are equivalent:
\begin{itemize}
\item[(i)] There exists a partition $\pi_{0}\in L $ such that $\pi_{0}\neq [n]$ and for every $\pi\in L $ we have that $\mu(\pi)=\mu(\pi\wedge \pi_{0})$,
\item[(ii)] $\mu_I=\mu(\pi_0(I))$ for all $I\subseteq [n]$,
\item[(iii)] $\ell(\pi)=0$ for all $\pi\not\leq\pi_0$,
\item[(iv)] $\ell_I=0$ unless $I$ is contained in a single block of $\pi_0$.
\end{itemize}
\end{thm}
\begin{proof}
The equivalence of (i) and (ii) follows from the fact that $\mu(\pi)$ is a multiplicative function of $L$. Hence (i)$\Rightarrow$(ii) follows by taking $\pi=[n]$ and then constraining to elements of $I$. The opposite implication follows by taking $I$ to be blocks of $\pi$. We now prove \mbox{(i)$\Rightarrow$(iii)}. Using Definition~\ref{def:xicumulants}  we obtain
\begin{eqnarray*}
 \ell(\nu)\,\,\,=\,\,\,\sum_{\pi\leq \nu}\mm(\pi,\nu)\mu(\pi)=\sum_{\pi\leq \nu}\mm(\pi,\nu)\mu(\pi\wedge \pi_{0})=\\
=\sum_{\delta\leq \nu\wedge \pi_0}\left(\sum_{\pi\wedge\pi_{0}=\delta}\mm(\pi,\nu)\right)\mu(\delta),
\end{eqnarray*}
where the inner sum in the last expression is over all $\pi$ in $[\hat{0},\nu]$  such that $\pi\wedge\pi_0=\delta$ (or $\pi\wedge(\pi_0\wedge\nu)=\delta$). To show (iii), we are interested only in $\nu\not\leq\pi_0$ and hence  we can assume that $\nu\neq\hat{0}$. The interval $[\hat{0},\nu]\subseteq L $ is a lattice with at least two elements, and, whenever $\nu\not\leq\pi_0$, also $\pi_{0}\wedge \nu\neq \nu$. Therefore, by Lemma~\ref{lem:sumzeroM} for all $\delta\leq\nu\wedge\pi_0$ the sum $\sum_{\pi\wedge\pi_{0}=\delta}\mm(\pi,\nu)$ vanishes. Hence $\ell(\nu)=0$ unless $\nu\leq\pi_0$.

To show \mbox{(iii)$\Rightarrow$(i)} note that if $\ell(\delta)=0$ for all $\delta\not\leq \pi_0$ then for every $\pi\in L$
$$
\mu(\pi)\,\,=\,\,\sum_{\delta\leq \pi}\ell(\delta)\,\,=\,\,\sum_{\delta\leq \pi\wedge \pi_0}\ell(\delta)\,\,=\,\,\mu(\pi\wedge \pi_0).
$$

To see that (iv) follows from  (i) and (iii), apply (i) with $L(I)$ in place of $L$. If $I$ is not contained in a block of $\pi_0$ then $\pi_0(I)$ is not the maximal element of $L(I)$  and by (i) this gives $\mu(\pi)=\mu(\pi\wedge \pi_0(I))$ for every $\pi\in L(I)$. Now $\ell_I=0$ by (iii). 

Finally we show that (iv) implies (ii) using induction with respect to $|I|$. If $I=\{i,j\}$ such that $i$ and $j$ lie in different blocks of $\pi_0$ then $\pi_0(I)=i|j$. Since $\ell_{ij}=\mu_{ij}-\mu_i\mu_j=0$, (ii) holds if $|I|=2$. Suppose now that (ii) holds for all $|I|<d$ and let now $I\subseteq [n]$ be such that $|I|=d$ and $\pi_0(I)\neq I$ (otherwise (ii) holds trivially). By  (\ref{xi-cumul}) we have
$$
\ell_I\,\,=\,\,\sum_{\pi\in L(I)}\mm(\pi)\mu(\pi).
$$
If $\pi<I$ then $\mu(\pi)$ is a product of some $\mu_B$, where $|B|<d$ and hence by assumption $\mu(\pi)=\mu(\pi\wedge \pi_0(I))$. We can rewrite the above equation as
\begin{equation}\label{eq:lAaux}
\ell_I\quad =\quad \mu_I-\mu(\pi_0(I))+\sum_{\pi\in L(I)}\mm(\pi)\mu(\pi\wedge \pi_0(I)).
\end{equation}
The last summand can be rewritten as 
$$\sum_{\delta\leq \pi_0(I)}\Big[\sum_{\pi\wedge\pi_0(I)=\delta}\mm(\pi)\Big]\mu(\delta),$$
which is zero by Lemma \ref{lem:weisner} because $\pi_0(I)\neq I$. Therefore, (\ref{eq:lAaux}) becomes $\ell_I=\mu_I-\mu(\pi_0(I))$. Since $\ell_I=0$ by assumption, we obtain that (ii) holds for $|I|=d$ and hence it holds for all $I\subseteq[n]$. 
\hfill$\Box$\end{proof}

This result gives an immediate corollary which generalizes property (P1) of the classical cumulants.
\begin{prop}\label{prop:P1}
Suppose there exists a partition $\pi_{0}\in L $ such that $\indep_{B\in \pi_0}X_B$. Then $\ell(\pi)=0$ for all $\pi\not\leq\pi_0$ or equivalently $\ell_A=0$ unless all the elements of $A$ are contained in a single block of $\pi_0$. \end{prop}

This proposition shows one of the important features of $L$-cumulants. For cumulants, by (P1),  all marginal independencies imply that $k_{1\cdots n}=0$. In the case of $L$-cumulants only some of the independencies imply vanishing (see Example \ref{ex:n3interv}). Hence, this new coordinate system can be designed to better fit the model under consideration. This concept will be explained in more detail for tree cumulants in Section \ref{sec:treecum}.

We formulate an additional condition on the family of lattices $\mathbf{L}$, which we require to hold only when this is explicitly stated.
\begin{itemize}
\item[(C1)] For  every $A\in \cA(\cX)$ and every $i\in A$ the split $i|(A\setminus i)$ is in $L(A)$.
\end{itemize}
Among the partitions in Definition~\ref{def:partitions} only the lattice of interval partitions does not satisfy (C1). 
\begin{prop}[Semi-invariance]\label{prop:translation}
Let $\mathbf{L}$ satisfy (C1) and $\widetilde{X}=X+a$, where $a\in\R^{n}$ is any constant vector. Denote by $\widetilde{\ell}_A$ the $L$-cumulant of $\widetilde{X}_A$. Then $\widetilde{\ell}_i={\ell}_i+a_{i}$ for all $i=1,\ldots, n$ and  $\widetilde{\ell}_A= {\ell}_A$ for any multiset $A\in \cA(\cX)$ such that $|A|\geq 2$. 
\end{prop}
\begin{proof}
Without loss of generality assume $A=[n]$. Since $a=\sum a_{i}e_{i}$, where the $e_{i}$'s are the unit vectors in $\R^{n}$, it suffices to prove this result only in the case when $a$ is such that $a_{1}$ is the only non-zero entry. In this case write $\widetilde{X}_{1}=X_{1}+a_{1}$ as $X_{1}-\mu_{1}+(a_{1}+\mu_{1})$, where $\mu_{1}=\E X_{1}$ and $a_{1}+\mu_{1}=\E \widetilde{X}_{1}$. Hence, if the split $\pi_{0}=1|\{2,\ldots n\} \,\in L $ then for every $\pi\in L$,
$$
\widetilde{\mu}(\pi)\,\,\,=\,\,\,\mu(\pi)-\mu(\pi\wedge \pi_{0})+\widetilde{\mu}(\pi\wedge \pi_{0}).
$$
 It follows that
\begin{equation}\label{eq:xi-prime-trans}
\begin{array}{l}
 \widetilde{\ell}_{1\cdots n}\quad=\quad\sum_{\pi\in  L }\mm(\pi)\mu(\pi)-\\
 \quad-\sum_{\pi\in  L }\mm(\pi)\mu(\pi\wedge \pi_{0})+\sum_{\pi\in  L }\mm(\pi)\widetilde{\mu}(\pi\wedge \pi_{0}).
\end{array} 
\end{equation}
Since $  L $ is a lattice and $\pi_{0}\neq [n]$, by Lemma~\ref{lem:sumzeroM} we have that  $\sum_{\pi\wedge\pi_{0}=\nu}\mm(\pi)=0$ for each $\nu\in L $ and hence the second and third summand in (\ref{eq:xi-prime-trans}) are zero. The proof is completed because the first summand is exactly ${\ell}_{1\cdots n}$. 
\hfill$\Box$\end{proof}

The following result shows that the central moments are $L$-cumulants induced by the lattice of one-cluster partitions $\cC([n])$.
\begin{prop}\label{prop:centralm}
Let $X$ be a random vector with values in $\cX$. Then the central moments $\mu_A'$ for $|A|\geq 2$ are equal to the corresponding $L$-cumulants induced by $\mathbf{C}=(\cC(A))_{A\in \cA(\cX)}$.
\end{prop}
\begin{proof}
Denote by $\mathfrak{c}$ the $L$-cumulants induced by the family $\mathbf{C}$ of one-cluster partition lattices. Let $A\in \cA(\cX)$ be such that $|A|\geq 2$. Since every split of the form $i|(A\setminus i)$ is a one-cluster partition, by Proposition \ref{prop:translation}, we can write $\mathfrak{c}$ in terms of the central moments
$$
\mathfrak{c}_{A}\quad =\quad \sum_{\pi\in \cC(A)}\mm(\pi)\prod_{B\in \pi}\mu'(B)\qquad\mbox{for all }|A|\geq 2.
$$
However, $\mu'_i=0$ for every $i\in [n]$ and hence the only non-zero term of the above sum is where $\pi=A$, which proves that $\mathfrak{c}_A=\mu'_A$.
\hfill$\Box$\end{proof}

The correspondence between the lattice of one-cluster partitions and central moments gives also the following explicit, simple and computationally efficient formula for central moments in terms of moments.
\begin{lem}\label{lem:centrmomform}
Let $X$ be a random vector with values in $\cX$. For every $A\in \cA(\cX)$ such that $|A|\geq 2$ we have:
\begin{equation}\label{eq:centralnice}
\mu_A'\quad=\quad\sum_{B\subseteq A}(-1)^{|A\setminus B|}\mu_B\prod_{i\in A\setminus B}\mu_i.
\end{equation}
\end{lem}
\begin{proof}
Use (\ref{eq:Mobiusfo1c}) and Proposition \ref{prop:centralm} to write 
$$
\mu'_A=\sum_{\hat{0}<\pi\in \cC(A)}(-1)^{|\pi|-1}\prod_{B\in \pi}\mu_B+(-1)^{|A|-1}(|A|-1)\prod_{i\in A}\mu_i.
$$
Let $B_0$ be the distinguished non-singleton block in each of the product $\prod_{B\in \pi}\mu_B$ above. Then $|\pi|-1=|A\setminus B_0|$. Hence, every $\prod_{B\in \pi}\mu_B$ corresponds to some $\mu_{B_0}\prod_{i\in A\setminus B_0}\mu_i$ in (\ref{eq:centralnice}) with the same coefficient. The remaining part is to check that the coefficient of $\prod_{i\in A}\mu_i$ is also the same, but this is an easy check.
\hfill$\Box$\end{proof}
 
\begin{exmp}
Let $A=\{1,1,2,2\}$ and list all multisubsets of $A$ as defined in the beginning of Section \ref{sec:XCdef}. We easily check that 
$$
\mu_{1122}'=\mu_{1122}-2\mu_1\mu_{122}-2\mu_2\mu_{112}+\mu_{11}\mu_2^2+4\mu_{12}\mu_1\mu_2+\mu_1^2\mu_{22}-3\mu_1^2\mu_2^2,
$$
which can be verified also by hand.
\end{exmp}

\subsection{Multilinear transformations}

By property (P3) cumulants behave nicely under multilinear transformations. In this section, to study similar properties for general $L$-cumulants, we restrict to $\mathbf{L}$ satisfying the following condition.
\begin{itemize}
\item[(C2)] For every $A\in \cA(\cX)$ the lattice $L(A)$ is isomorphic to $L{([d])}$, where $d=|A|$.
\end{itemize}
This property is satisfied by construction for  $\mathbf{\Pi}$, $\mathbf{I}$, $\mathbf{NC}$, and $\mathbf{C}$. If (C2) holds then, for every $d$-tuple $(i_1,\ldots,i_d)\in [n]^d$ we define $\cL^{(d)}$ as a $n\times\cdots\times n$ tensor of the form
\begin{equation}\label{eq:generalLcumul}
\cL^{(d)}_{i_1\cdots i_d}\,\,\,=\,\,\,\sum_{\pi\in L([d])}\mm(\pi)\prod_{B\in \pi}\mu_{i_B}.
\end{equation}
Note that in general $\cL^{(d)}_{i_1\cdots i_d}$ may differ from $\ell_{i_1\cdots i_d}$. For example if $L=\cI([3])$ then $\ell_{213}=\ell_{123}$ because the definition of $L$-cumulants does not depend on the ordering of the elements in $[n]$. On the other hand, we have $\cL^{(d)}_{123}\neq \cL^{(d)}_{213}$ because
$$\cL^{(d)}_{123}\,\,=\,\,\mu_{123}-\mu_1\mu_{23}-\mu_{12}\mu_3+\mu_1\mu_2\mu_3$$
and
$$
\cL^{(d)}_{213}\,\,=\,\,\mu_{123}-\mu_2\mu_{13}-\mu_{12}\mu_3+\mu_1\mu_2\mu_3.
$$ 

  The following proposition shows that the tensor $\cL^{(d)}$, for any $d\geq 1$, under linear mappings transforms as a contravariant tensor.  
\begin{prop}\label{prop:P3}
Let $X=(X_{1},\ldots, X_{n})$ be a random vector. Consider $L$-cumulants defined by $\mathbf{L}$ satisfying (C2). Let  $Q=[q_{ij}]\in \R^{m\times n}$ and $\widetilde{X}=QX\in \R^m$. Define $[\widetilde{\mu}_A]$, $[\widetilde{\ell}_A]$ and $\widetilde{\mathcal{L}}^{(d)}$ as counterparts of $[\mu_A]$, $[\ell_A]$ and $\mathcal{L}^{(d)}$ for $\widetilde{X}$ accordingly. Then for each $d\geq 1$, $$\widetilde{\mathcal{L}}^{(d)} \quad=\quad Q\cdot {\mathcal{L}}^{(d)},$$ where $Q\cdot {\mathcal{L}}^{(d)}$ is the multilinear action on a $d$-dimensional tensor defined by 
\begin{equation}\label{eq:multilinL}
(Q\cdot  {\mathcal{L}}^{(d)})_{i_{1}\cdots i_{d}}\quad =\quad \sum_{j_{1}=1}^{n}\cdots\sum_{j_{d}=1}^{n} q_{i_{1}j_{1}}\cdots q_{i_{d}j_{d}}  \cL^{(d)}_{j_{1}\cdots j_{d}}
\end{equation}
for each $d\geq 1$ and $i_{1},\ldots, i_{d}\in [m]$. 
\end{prop}
\begin{proof}
By (\ref{eq:generalLcumul}) we have
$$
(Q\cdot {\mathcal{L}}^{(d)})_{i_{1}\cdots i_{d}}=\sum_{j_{1}=1}^{n}\cdots\sum_{j_{d}=1}^{n} q_{i_{1}j_{1}}\cdots q_{i_{d}j_{d}}\left( \sum_{\pi\in  L([{d}])}\mm(\pi)\prod_{B\in \pi}\mu_{j_B}\right).
$$
Write $\mu_{j_B}$ explicitly as $\E\left[\prod_{b\in B} X_{j_{b}}\right]$. Then, using (C2), after changing the ordering  of products and summations we obtain
$$
(Q\cdot {\mathcal{L}}^{(d)})_{i_{1}\cdots i_{d}}=\sum_{\pi\in  L([{d}])}\mm(\pi)\prod_{B\in \pi} \E \left[\prod_{b\in B}(\sum_{j_{b}=1}^{n}q_{i_{b}j_{b}}X_{j_{b}}) \right].
$$
Since $\widetilde{X}_{i_{b}}=\sum_{j_{b}=1}^{n}q_{i_{b}j_{b}}X_{j_{b}}$ we obtain
$$
(Q\cdot {\mathcal{L}}^{(d)})_{i_{1}\cdots i_{d}}=\sum_{\pi\in  L([{d}])}\mm(\pi)\prod_{B\in \pi} \widetilde{\mu}_{i_B}=( \widetilde{\mathcal{L}}^{(d)})_{i_{1}\cdots i_{d}},
$$
which finishes the proof.
\hfill$\Box$\end{proof}
Although for some $\mathbf{L}$ the property (P3) may not hold, the homogeneity holds for all $L$-cumulants. Thus, if $\widetilde{X}=(\lambda_{1}X_{1},\ldots, \lambda_{n}X_{n})$ for some $\lambda=(\lambda_{1},\ldots,\lambda_{n})\in (\R^{*})^{n}$   then $\widetilde{\ell}_{A}=\prod_{i\in A}\lambda_{i} \ell_{A}$ for every $A\in \cA(\cX)$.

\subsection{Conditional $ L$-cumulants}\label{sec:XCconditional} Suppose we are given the conditional cumulants of $X=(X_{1},\ldots, X_{n})$ conditional on some random variable $Y$ and we want to obtain the unconditional cumulants. This is a common problem with hidden variable models. On the level of moments this relationship is straightforward since
$$\mu_{A}\,\,\,=\,\,\,\E \big[\prod_{i\in A}X_{i}\big]\,\,\,=\,\,\,\E\Big[ \E\big[\prod_{i\in A}X_{i}|Y\big]\Big]$$ for every multiset $A$. For cumulants, or more generally for $ L$-cumulants, the situation is a bit more complicated. 

For every multiset $A\in \cA(\cX)$ denote by $k_A^Y$ the \textit{conditional cumulant of $X_A$ given $Y$}, that is a cumulant computed as in Definition \ref{def:xicumulants} but with moments replaced by conditional moments. Note that each $k_A^Y$ is itself a random variable. For any $\pi\in \Pi(A)$, by $\hat{k}_\pi$ denote the cumulant of the random vector $(k^Y_B)_{B\in \pi}$. It is  known from  \cite{brillinger1969calculation} that for every $A\in \cA(\cX)$:
\begin{equation}\label{eq:cond-cumul0}
 {k}_{A}\quad=\quad\sum_{\pi\in  L(A)} \hat{k}_\pi.
\end{equation}
This in particular generalizes the well-known formula
$$
{\rm Cov}(X,Z)\quad =\quad \E[{\rm Cov}(X,Z|Y)]+{\rm Cov}(\E[X|Y],\E[Z|Y]).
$$
In Theorem~\ref{th:conditional} we give a purely combinatorial proof of (\ref{eq:cond-cumul0}). For our purposes it is slightly more constructive than a similar proof of the same result in \cite{speed1983cumulants}. Also it  immediately enables us to formulate this result  for $ L$-cumulants in the case when $\mathbf{L}$ satisfies the following property.
\begin{itemize}
\item[(C3)] For every $n\geq 0$ and each $\pi\in L $ the interval $[\pi,[n]]\subseteq L $ is isomorphic to $L{([|\pi|])}$. 
\end{itemize}
This property is satisfied for $\mathbf{\Pi}$ (see \cite[Example~3.10.4]{stanley2006enumerative}). A sufficient condition for $\mathbf{L}$ to satisfy (C3) is that for every $n\geq 0$ the lattice $L $ forms a join subsemilattice of $\Pi([n])$. Therefore, $\textbf{I}$ as well as the lattice of tree partitions for sufficiently regular trees (for example caterpillars) both satisfy the property. Condition (C3) does not hold however for the lattice of one-cluster partitions, (general) tree partitions and non-crossing partitions.  
 
 For every multiset $A\in \cA(\cX)$ denote by $\ell_A^Y$ the conditional $L$-cumulant of $X_A$ given $Y$. For any $\pi\in L(A)$, by $\hat{\ell}_\pi$ denote the $L$-cumulant of the random vector $(\ell^Y_B)_{B\in \pi}$.
\begin{thm}[Brillinger's formula for $L$-cumulants]\label{th:conditional}
Let $X=(X_{1},\ldots, X_{n})$ be a random vector and $Y$ be a random variable. If $\mathbf{L}$ satisfies  (C3) then 
\begin{equation*} \ell_{1\cdots n}\quad=\quad\sum_{\pi\in  L } \hat{\ell}_\pi.
\end{equation*}
\end{thm}
\begin{proof}
Since $\mu_{B}=\E\mu_{B}^{Y}$, by (\ref{eq:inversexi}) we obtain the identity
\begin{equation}\label{eq:cond-cum-ident1}
\mu_{B}\,\,\,=\,\,\,\E\mu_{B}^{Y}\,\,\,=\,\,\,\sum_{\delta\in L(B)}\E\big[\prod_{C\in\delta} \ell_{C}^{Y}\big].
\end{equation}
Using (\ref{xi-cumul}) and replacing (\ref{eq:cond-cum-ident1}) for each $\mu_{B}$ we can write
$$
\begin{array}{rcl}
 \ell_{1\cdots n}&=&\sum_{\pi\in L }\mm(\pi)\prod_{B\in\pi} \left( \sum_{\delta\in L(B)}\E\left[\prod_{C\in\delta} \ell_{C}^{Y}\right]\right)=\\
&=& \sum_{\pi\in L }\mm(\pi)\sum_{\delta\leq \pi}\prod_{B\in\pi} \E\left[\prod_{C\in\delta(B)} \ell_{C}^{Y}\right],
\end{array}
$$
where $\delta(B)$ denotes the partition $\delta\in L $ constrained to $B\in \pi$. We change the order of summation to obtain
\begin{equation}\label{eq:cond-cum2}
 \ell_{1\cdots n}\,\,\,=\,\,\,\sum_{\delta\in L }\left(\sum_{\pi\geq \delta}\mm(\pi)\prod_{B\in\pi} \E\big[\prod_{C\in\delta(B)} \ell_{C}^{Y}\big]\right).
\end{equation}
For each $\delta=C_{1}|\cdots|C_{r}\in L $ denote the set of its blocks by $\mathbb{B}_{\delta}=\{C_{1},\ldots, C_{r}\}$. By (C3) the interval $[\delta,[n]]$ is isomorphic to $L(\mathbb{B}_{\delta})$ which is isomorphic to $L([{|\delta|}])$ and hence the expression in brackets in (\ref{eq:cond-cum2}) can be rewritten as 
$$
\sum_{\nu\in L({\mathbb{B}_{\delta}})}\mm_{\mathbb{B}_{\delta}}(\nu,{\mathbb{B}_{\delta}})\prod_{B\in \nu} \E\big[\prod_{C\in B}  \ell_{C}^{Y}\big],
$$
which by definition is equal to $\hat{\ell}_\delta$. 
\hfill$\Box$\end{proof}


If (C3) does not hold  and we want to perform some efficient conditional computations, we can still use the classical Brillinger's formula for cumulants and then translate them back to $L$-cumulants using Proposition~\ref{prop:LinK}. Moreover, for some special statistical models the following result may be useful. It works for all families $\mathbf{L}$.
\begin{prop}\label{prop:indep}
Let $X=(X_{1},\ldots, X_{n})$ be a random vector and $Y$ a random variable. If $X_{1}\indep\cdots\indep X_{n}|Y$, then 
\begin{equation*}
 \ell_{1\ldots n}\quad=\quad\hat{\ell}_{1|2|\cdots|n},
\end{equation*}
where by definition $\hat{\ell}_{1|2|\cdots|n}$  is the $L$-cumulant of the random vector $(\ell_1^Y,\ldots,\ell^Y_n)=(\mu_1^Y,\ldots,\mu^Y_n)$.
\end{prop}
 \begin{proof}
Since $X_{1}\indep \cdots\indep X_{n}|Y$, by Proposition~\ref{prop:P1}, $\ell_{C}^{Y}=0$ unless $|C|=1$. Moreover, we have 
\begin{equation*}
\mu_{B}\,\,\,=\,\,\,\E[\prod_{i\in B}\mu_{i}^{Y}].
\end{equation*}
 Using (\ref{xi-cumul}) and replacing the above identity for each $\mu_{B}$ we can write
$$
\begin{array}{rcl}
 \ell_{1\cdots n}&=&\sum_{\pi\in L }\mm(\pi)\prod_{B\in\pi} \E\left[\prod_{i\in B} \mu_{i}^{Y}\right].
\end{array}
$$
But since $\ell^Y_i=\mu^Y_i$, the right hand side in the above equation is exactly the $L$-cumulant of the random vector $(\ell_1^Y,\ldots,\ell^Y_n)$.
\hfill$\Box$\end{proof}

To see how this result may be relevant in geometry see  Example~\ref{ex:secant}.

\section{Tree cumulants and hidden Markov processes}\label{sec:treecum}

In this section we complement the discussion of tree cumulants and show how they can be used to analyze more general processes on trees. 

\subsection{Tree models}
Let $T^{r}$ be a rooted tree with vertex set $V$ and edge set $E$, that is  a tree with one distinguished node $r\in V$ called the root and all the edges directed away from $r$. Let $X=(X_{v})_{v\in V}$  be a vector of  binary random variables with values $0$ and $1$. Consider the Bayesian network for $X$ represented by $T^{r}$. Each node $v$ corresponds to a random variable $X_v$ and the structure of $T^r$ imposes some constraints on the joint distribution of $X$ (see for example \cite{lauritzen:96}). Define $\cM_{T}$ as the model obtained from this Bayesian network by taking the marginal distributions over the leaves of $T^{r}$. We call $\cM_{T}$ the \textit{two-state general Markov model} (for example  \cite[Chapter 8]{semple2003pol}). We omit the rooting in the notation because the model does not depend on the rooting. In other words, for any alternative rooting the induced parametrization will lead to the same model. 

The parametric formulation of the model is obtained by expressing the marginal distribution of $X$ over the leaves of $T^r$ in terms of the marginal distribution of the root $r$ and conditional distributions of each $v\in V\setminus \{r\}$ given its parent in $T^r$ denoted by ${\rm pa}(v)$. Assume that $T^r$ has $n$ leaves and label  them by elements of $[n]$. The distribution over the set of leaves satisfies 
\begin{equation}\label{eq:gmm-param}
p({x_1,\ldots, x_n})\quad=\quad\sum_\cH p_r(x_r)\prod_{v\in V\setminus r} p_{v|{\rm pa}(v)}(x_v|x_{{\rm pa}(v)}),
\end{equation}
where $\cH$ is the set of all $x\in\{0,1\}^{V}$ such that the restriction to the leaves of $T$ is equal to $(x_1,\ldots,x_n)$. The model is given as the image of (\ref{eq:gmm-param}) in $\Delta_\cX$, where each point corresponds to a different choice of values for conditional probabilities on the right hand side of this parametrization. If $m$ denotes the number of inner nodes of $T$ then this parametrization has $2^{m}$ terms. For large trees this is a big polynomial which complicates the geometric and algebraic analysis of these models. 

The two-state general Markov model can be equivalently defined by a set of conditional independence statements. This follows from the general theory of graphical models (see \cite[Section 3.2.2]{lauritzen:96}). We say that two disjoint subsets $A,B$ of the set of vertices $V$ of $T$ are \textit{separated} by another subset $C$ if every undirected path from a node in $A$ to a node in $B$ necessarily crosses $C$. The set of all conditional independence statements which define the general Markov model are given by all $A\indep B|C$ for all disjoint subsets $A,B,C\subseteq V$ such that $C$ separates $A$ and $B$. For example the $4$-star tree model discussed in Section \ref{sec:tressforsecs} is defined by $X_1\indep X_2\indep X_3\indep X_4|Y$ because the inner node separates all the leaves from each other. 

Before we recall the main result of \cite{pwz-2010-identifiability}, let us give some intuition on why tree cumulants may be helpful in the study of tree models. Suppose that for some edge $(u,v)$ in $T^r$ we impose on the model $\cM_T$ that in addition $X_u\indep X_v$. This corresponds to removing the edge $(u,v)$ from $T^r$ and considering the model of the induced forest. Let $A|B$ be the split of the set of leaves $[n]$ induced by removing the edge $(u,v)$. Then the independence statement $X_u\indep X_v$ implies also that $X_A\indep X_B$.
\begin{exmp}\label{ex:aboveexample}
Let $T$ be the quartet tree in Figure  \ref{fig:quartet} rooted in $a$. The independence $(X_1,X_2)\indep (X_3,X_4)$ defines a valid submodel of the tree model for $T$. This submodel is defined by requesting $X_a\indep X_b$ and hence it is given as the image of the subspace of the parameter space restricted to $p_{b|a}({1|0})=p_{b|a}({1|1})$. 
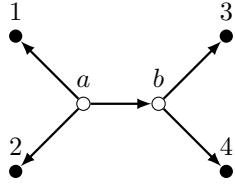
\begin{figure}[h!]
\centering
\tikzstyle{vertex}=[circle,fill=black,minimum size=5pt,inner sep=0pt]
\tikzstyle{hidden}=[circle,draw,minimum size=5pt,inner sep=0pt]
  \begin{tikzpicture}
  \node[vertex] (1) at (-.9,.9)  [label=above:$1$] {};
    \node[vertex] (2) at (-.9,-.9)  [label=above:$2$] {};
    \node[vertex] (3) at (1.9,.9) [label=above:$3$]{};
    \node[vertex] (4) at (1.9,-.9) [label=above:$4$]{};
    \node[hidden] (a) at (0,0) [label=above:$a$] {};
    \node[hidden] (b) at (1,0) [label=above:$b$] {};
          \draw[->,-latex,line width=.3mm] (a) to (b);
    \draw[->,-latex,line width=.3mm] (a) to (1);
    \draw[->,-latex,line width=.3mm] (a) to (2);
    \draw[->,-latex,line width=.3mm] (b) to (3);
    \draw[->,-latex,line width=.3mm] (b) to (4);
  \end{tikzpicture}
  \caption{A quartet tree.}\label{fig:quartet}
\end{figure}
\end{exmp}

By Proposition \ref{prop:binindep} there exists a tree partition $\pi_0$ such that $\indep_{B\in \pi_0}X_B$ if and only if $\mathfrak{t}_I=0$ whenever $I\subseteq[n]$ is not completely contained in one of the blocks of $\pi_0$. In Example \ref{ex:aboveexample}, because $12|34$ is a valid tree partition, the marginal independence $(X_1,X_2)\indep (X_3,X_4)$ holds if and only if $\mathfrak{t}_I=0$ for all $I\subseteq \{1,2,3,4\}$ such that $I$ is not contained neither in $\{1,2\}$ nor $\{3,4\}$. Hence all $\mathfrak{t}_{13}$, $\mathfrak{t}_{14}$, $\mathfrak{t}_{23}$, $\mathfrak{t}_{24}$, $\mathfrak{t}_{134}$,$\mathfrak{t}_{234}$,$\mathfrak{t}_{123}$,$\mathfrak{t}_{124}$ and $\mathfrak{t}_{1234}$ vanish whenever $p_{b|a}({1|0})=p_{b|a}({1|1})$. These kind of considerations help to understand why tree cumulants are helpful for describing the two-state general Markov models. They also help to intuitively understand the result in Theorem \ref{th:maintree}, which we now state formally. 

Let $\eta_{uv}=p_{v|u}({1|1})-p_{v|u}({1|0})$ for each $(u,v)\in E$. As we have shown $\eta_{uv}=0$ if and only if $X_u\indep X_v$. Moreover, let $\bar{\mu}_{v}=1-2\mu_v$ for $v\in V$. 
\begin{thm}[Zwiernik, Smith \cite{pwz-2010-identifiability}]\label{th:maintree}
Let $T$ be trivalent tree. Then the two-state general Markov model $\cM_{T}$ can be equivalently expressed in the space of tree cumulants by $\ell_{i}=\mu_{i}=\frac{1}{2}(1-\bar{\mu}_{i})$ for $i=1,\ldots, n$; and for all $|I|\geq 2$
\begin{equation*}
\mathfrak{t}_{I}\quad=\quad\frac{1}{4}(1-\bar{\mu}_{r(I)}^{2})\prod_{\deg(v)=3}\bar{\mu}_{v}\prod_{(u,v)\in E(I)}\eta_{uv},
\end{equation*}
where $V(I)$ and $E(I)$ denotes vertex and edge sets of the tree $T(I)$, $r(I)$ is the root of $T(I)$ and $\deg(v)$ denotes the valency of $v$ in $T(I)$.\end{thm}
\begin{exmp}\label{ex:quartet}
Let $T^r$ be a quartet tree in Figure \ref{fig:quartet}. Then by Theorem~\ref{th:maintree} we have for example $\mathfrak{t}_{12}=\frac{1}{4}(1-\bar{\mu}_{a}^{2})\eta_{a1}\eta_{a2}$, $\mathfrak{t}_{13}=\frac{1}{4}(1-\bar{\mu}_{a}^{2})\eta_{a1}\eta_{ab}\eta_{b3}$, $\mathfrak{t}_{34}=\frac{1}{4}(1-\bar{\mu}_{b}^{2})\eta_{b3}\eta_{b4}$ and
$$
\mathfrak{t}_{1234}\,\,\,=\,\,\,\frac{1}{4}(1-\bar{\mu}_{a}^{2})\bar{\mu}_{a}\bar{\mu}_{b}\eta_{1a}\eta_{2a}\eta_{ab}\eta_{b3}\eta_{b4}.
$$
We also infer from this that $$\mathfrak{t}_{I\cup J}\mathfrak{t}_{I'\cup J'}-\mathfrak{t}_{I\cup J'}\mathfrak{t}_{I'\cup J}\,\,\,=\,\,\,0$$ for all $I,I'\in \{\{1\},\{2\},\{1,2\}\}$ and $J,J'\in \{\{3\},\{4\},\{3,4\}\}$.
\end{exmp}

Theorem \ref{th:maintree} can be applied only for trivalent trees and hence it does not hold for $n$-star tree models discussed earlier (see also Remark \ref{rem:naively}). We can use however the fact that any non-trivalent tree model is a submodel of some model of a trivalent tree. Thus, if $T$ is not trivalent then we take any trivalent tree $T^*$ such that $T$ can be obtained from $T^*$ by edge contractions. Now the two-state general Markov model for $T$, when expressed in tree cumulants of $T^*$, is parametrized by $\ell_i=\mu_i$ for $i=1,\ldots,n$, and for all $|I|\geq 2$
 \begin{equation*}
 \mathfrak{t}_{I}\quad=\quad\frac{1}{4}(1-\bar{\mu}_{r(I)}^{2})\prod_{v\in V(I)\setminus I}\bar{\mu}_{v}^{\deg(v)-2}\prod_{(u,v)\in E(I)}\eta_{uv}.
\end{equation*}
   In the quartet tree of Example~\ref{ex:quartet} we can contract the edge $(a,b)$ to obtain the $4$-star tree in Figure \ref{fig:4star}. This contraction corresponds to the subspace of the parameter space given by $\bar{\mu}_a=\bar{\mu}_b$ and $\eta_{ab}=1$. This induces the parametrization of the secant variety given in Example~\ref{ex:secant}. The same can be obtained for any $n$-star tree model with $n\geq 4$. For more details see \cite{pwz-2010-identifiability}.


\subsection{Binary hidden Markov processes}

We now show that tree cumulants can be useful also for other related statistical models. We consider models with an underlying two-state Markov chain which is not observed, where the observed variables are independent given this Markov chain. An example is given by the hidden Markov model or some simple cases of Markov switching models without autoregressive terms (see for example \cite{hamilton1989new}). In this section we refer to all these models as \textit{binary hidden Markov processes}. 

Consider tree cumulants induced by the caterpillar tree and define the \textit{normalized tree cumulants} as 
$$
\bar{\mathfrak{t}}_{I}\,\,\,=\,\,\,\prod_{i\in I}\frac{1}{\sqrt{k_{ii}}} \mathfrak{t}_{I}\qquad\mbox{for all } I\subseteq [n],
$$
which is always well defined if all the variables in the system are non-degenerate (a degenerate random variable takes only one value with nonzero probability). With this definition $\rho_{ij}:=\bar{\mathfrak{t}}_{ij}$ is just the usual \textit{correlation} between $X_{i}$ and $X_{j}$, and $\gamma_i:=\bar{\mathfrak{t}}_{iii}$ is the \textit{skewness} of $X_i$.

In this section we deal with an observed vector $X=(X_1,\ldots, X_n)$ and a hidden vector $H=(H_1,\ldots,H_n)$. Since we need to consider mixed tree cumulants involving indices from both vectors, we introduce the following convention. Whenever an index involves $i$ referring to $H_i$ we write it as $\underline{i}$. Hence for example $k_{\underline{i}j}={\rm Cov}(H_i,X_j)$, $k_{\underline{i}\underline{i}}={\rm Var}(H_i)$, $k_{ii}={\rm Var}(X_i)$, $\gamma_{\underline{i}}=\E(H_i-\E H_i)^3/{\rm Var}(H_i)^{3/2}$, and $\mu'_{\underline{B}}=\E[ \prod_{i\in B}(H_i-\E H_i)]$.

It is well known that for every random variable $X$, if $Y$ is binary, then
\begin{equation}\label{eq:condexpbin}
\E(X|Y)\,\,\,=\,\,\,\E X+{\rm Cov}(X,Y)({\rm Var}(Y))^{-1}(Y-\E Y),
\end{equation}
where ${\rm Cov}(X,Y)({\rm Var}(Y))^{-1}$ is the linear regression coefficient of $X$ with respect to $Y$. The following proposition shows that the hidden Markov process has an elegant formulation  and all its normalized tree cumulants are parametrized by correlations and skewnesses. 
\begin{prop}\label{prop:HMPs}
Let $X=(X_{1},\ldots,X_{n})$ be a random vector and $H=(H_{1},\ldots,H_{n})$ a binary random vector (both non-degenerate). Assume that $X_{1}\indep\ldots \indep X_{n}|H$ and the conditional distribution of $X_{i}$ given $H$ depends only on $H_{i}$ for $i=1,\ldots,n$. Moreover, let $H$ form a Markov chain. Then for every $I=\{i_1,\ldots,i_d\}$ such that $1\leq i_1<\cdots<i_d\leq n$ the corresponding normalized tree cumulant satisfies
\begin{equation*}
\bar{\mathfrak{t}}_{I}\quad=\quad\prod_{j=2}^{d-1} \gamma_{\underline{i_j}}\prod_{i=i_1}^{i_d-1}\rho_{\underline{i}\,\underline{i+1}}\prod_{i\in I}\rho_{\underline{i}\, i}.
\end{equation*} 
\end{prop}
\begin{proof}
Before we prove the proposition we formulate the following result.
\begin{lem}\label{lem:abCind}
Suppose that $X=(X_1,\ldots,X_n)$ is a binary random vector such that $i\indep j\indep C|r$ for some disjoint $i,j,r\in [n]$ and $C\subseteq [n]$. Let $\eta_{rA}=\mu'_{rA}k_{rr}^{-1}$ for every $A\subseteq [n]$  and $\tau_{r}=k_{rrr}k_{rr}^{-1}$. Then
$$
\begin{array}{l}
\mu_{ijC}'\,\,\,=\,\,\,\eta_{ri}\eta_{rj}k_{rr}\mu_{C}'+\eta_{ri}\eta_{rj}\mu_{rC}'\tau_{r}.
\end{array}
$$
\end{lem}
\begin{proof}
Let $U_A:=\prod_{i\in A}(X_i-\E X_i)$ for every $A\subseteq [n]$. The conditional independence ${i}\indep {j}\indep {C}|{r}$ implies
$$
\mathbb{E}[U_{ijC}|U_{r}]\,\,\,=\,\,\,\mathbb{E}[U_{i}|U_{r}]\,\, \mathbb{E}[U_{j}|U_{r}] \,\,\mathbb{E}[U_{C}|U_{r}].
$$
Using (\ref{eq:condexpbin}) for the conditional expectations on the right hand side and then taking expectations on both sides yields
$$
\mu_{ijC}'\,\,\,=\,\,\,\eta_{ri}\eta_{rj}k_{rr}\mu'_{C}+\eta_{ri}\eta_{rj}\eta_{rC}k_{rrr}.
$$
Replace $\eta_{rC}=\mu_{rC}'k_{rr}^{-1}$ to obtain the formula in Lemma \ref{lem:abCind}.
\hfill$\Box$\end{proof}

To prove Proposition \ref{prop:HMPs} first assume that $I=[n]$ and by $L$ denote the lattice of tree partitions of the caterpillar tree with $n$ leaves. We can divide the partitions in $L$ into two groups:
\begin{itemize}
\item[1.] partitions with $1$ and $2$ in two different blocks $1A$ and $2B$, and
\item[2.] partitions with $1$ and $2$ in a single block $12A$
\end{itemize}

By Remark \ref{rem:centraltree}
 we can write
\begin{equation}\label{eq:auxtreecum}
\mathfrak{t}_{1\dots n}\,\,\,=\,\,\,\sum_{\pi\in L}\mathfrak{m}(\pi)\prod_{B\in \pi}\mu_{B}'.
\end{equation}
In the first group of partitions we always have either $A=\emptyset$ or  $B=\emptyset$. Since $\mu'_{1}=\mu'_{2}=0$, for every $\pi$ in the first group the corresponding summand in (\ref{eq:auxtreecum}) is zero. Let $\delta_0=12|3|\cdots|n$. The set of all partitions in the second group forms an interval $[\delta_0,[n]]$, which is isomorphic to the set of all tree partitions of the subtree $T_{2}$ of $T$ with $n-1$ leaves given by the hidden vertex $\underline{2}$ and the remaining leaves of $T$: $3,\ldots, n$. This isomorphism is given by replacing each block $12A$ with a block $\underline{2}A$. Denote the lattice of all partitions of $T_2$ by $L_2$. Since $[\delta_0,[n]]\simeq L_2$, the M\"{o}bius function on $L$ restricted to this interval is equal to the M\"{o}bius function on $L_2$. 
 
For every $A\subseteq [n]\setminus\{1,2\}$ we have that $X_{1}\indep X_{2}\indep X_{A}|H_2$ and hence, by Lemma \ref{lem:abCind}
$$
\mu_{12A}'\,\,\,=\,\,\,\eta_{\underline{2}1}\eta_{\underline{2}2}k_{\underline{2}\underline{2}}\mu_{A}'+\eta_{\underline{2}1}\eta_{\underline{2}2}\mu_{\underline{2}A}'\tau_{\underline{2}}.
$$
Therefore, (\ref{eq:auxtreecum}) becomes 
\begin{equation}\label{eq:treecumauxsum}
\mathfrak{t}_{1\dots n}\,\,\,=\,\,\,\sum_{\pi\in [\delta_0,[n]]}\mathfrak{m}(\pi)\prod_{B\in \pi}\mu_{B}'\cdot\eta_{\underline{2}1}\eta_{\underline{2}2}k_{\underline{2}\underline{2}}\mu_{A}'+\eta_{\underline{2}1}\eta_{\underline{2}2}\tau_{\underline{2}}\sum_{\pi\in L_2}\mathfrak{m}_{2}(\pi)\prod_{B\in \pi}\mu_{B}'.
\end{equation}
Let $\pi_{0}$ be a split $12|[n]\setminus\{1,2\}$. For every $\pi\in [\delta_0,[n]]$ the partition $\pi\wedge\pi_{0}$ is the partition obtained from $\pi$ by splitting the block $12A$ into two blocks $12$ and $A$. With this notation the first summand in (\ref{eq:treecumauxsum}) can be rewritten as
$$
\sum_{\nu\in[\delta_0,\pi_0]}\Big[\sum_{\pi:\pi\wedge\pi_{0}=\nu}\mathfrak{m}(\pi)\Big]\prod_{B\in \nu}\mu_{B}'\cdot\eta_{\underline{2}1}\eta_{\underline{2}2}k_{\underline{2}\underline{2}}.
$$
Since the interval $[\delta_0,[n]]$ forms a lattice then by Lemma \ref{lem:weisner} the above expression is zero. Since  $\sum_{\pi\in L_2}\mathfrak{m}(\pi)\prod_{B\in \pi}\mu_{B}'=\mathfrak{t}_{([n]\setminus\{1,2\})\cup \{a\}}$, then the second summand in (\ref{eq:treecumauxsum}) is
\begin{equation}\label{eq:auxtreeee}
\mathfrak{t}_{1\cdots n}\,\,\,=\,\,\,\eta_{\underline{2}1}\eta_{\underline{2}2}\tau_{\underline{2}}\mathfrak{t}_{([n]\setminus\{1,2\})\cup \{\underline{2}\}}.
\end{equation}
Using (\ref{eq:condexpbin}) we can also prove that  $\eta_{\underline{2}1}=k_{\underline{1}\underline{1}}\eta_{\underline{1}1}\eta_{\underline{1}\underline{2}}k_{\underline{2}\underline{2}}^{-1}$ (use the fact that $X_1\indep H_2|H_1$). In the next step we can apply the same procedure as above to express $\mathfrak{t}_{([n]\setminus\{1,2\})\cup \{\underline{2}\}}$ in (\ref{eq:auxtreeee}) in terms of $k_{\underline{2}\underline{2}}$, $\eta_{\underline{2}\underline{3}}$, $\eta_{\underline{3}{3}}$, $\tau_{\underline{3}}$ and $\mathfrak{t}_{([n]\setminus\{1,2,3\})\cup \{\underline{3}\}}$. We can do it recursively until we obtain
$$
\mathfrak{t}_{1\cdots n}=k_{\underline{1}\underline{1}}\prod_{i=1}^n\eta_{\underline{i}i}\prod_{i=1}^{n-1}\eta_{\underline{i}\,\underline{i+1}}\prod_{i=2}^{n-1}\tau_{\underline{i}}.
$$
Divide both sides by $\sqrt{k_{11}\cdots k_{nn}}$. The main proposition follows for $I=[n]$ after some obvious algebraic rearrangements. In the general case we first use the formula for $[n]$ to conclude that  for every $I=\{i_1,\ldots,i_d\}$, where $i_1<\ldots<i_d$
\begin{equation*}
\bar{\mathfrak{t}}_{I}\quad=\quad\prod_{j=2}^{d-1} \gamma_{\underline{i_j}}\prod_{j=1}^{d-1}\rho_{\underline{i_j}\,\underline{i_{j+1}}}\prod_{i\in I}\rho_{\underline{i}\, i}.
\end{equation*} 
To prove the final formula, it remains to show that 
$$\rho_{\underline{i_j}\,\underline{i_{j+1}}}=\prod_{i=i_j}^{i_{j+1}-1}\rho_{\underline{i}\,\underline{i+1}},$$
which can be proved by induction using (\ref{eq:condexpbin}).
\hfill$\Box$\end{proof}

This proposition  enables us to analyze the moment structure of hidden Markov processes. 
\begin{exmp}[Homogeneous binary hidden Markov model]
Consider a homogeneous binary hidden Markov model. In this case  ${H}=(H_{i})_{i=1}^{n}$ forms a homogeneous two-state Markov chain which we assume to  start from its stationary distribution. Moreover, the conditional distribution of $X_i$ given $H_i$ is the same for every $i=1,\ldots,n$. Under these assumptions the marginal distribution of $H_1$ is equal to the marginal distribution of $H_i$ for every $i=2,\ldots,n$.  Let $\gamma$ be the skewness of $H_1$, $\rho={\rm Corr}(H_{1},H_{2})$ be the one step correlation of the Markov chain ${H}$, and $b={\rm Corr}(H_{1},X_{1})$. By Proposition \ref{prop:HMPs}, for every $d\geq 2$ and $1\leq i_1<\ldots< i_d\leq n$,
\begin{equation}\label{eq:polellesx2}
\overline{\mathfrak{t}}_{i_{1}\cdots i_{d}}\quad=\quad b^{d}\rho^{i_{d}-i_{1}}\gamma^{d-2}.
\end{equation}
This in turn induces some constraints on the tree cumulants of the observed variables which may be useful to construct simple diagnostic tests for this class of models. For example it is easy to check that
$$
\overline{\mathfrak{t}}_{i (i+2)}\overline{\mathfrak{t}}_{j(j+2)}\,\,\,=\,\,\,\overline{\mathfrak{t}}_{k(k+3)}\overline{\mathfrak{t}}_{l(l+1)}\quad\mbox{for every } i,j,k,l=1,\ldots,n
$$
and that $\overline{\mathfrak{t}}_{ij}\overline{\mathfrak{t}}_{ik}\overline{\mathfrak{t}}_{jk}\geq 0$ for all $i<j<k$. The monomial parametrization in (\ref{eq:polellesx2}) enables us to obtain the equations for higher order tree cumulants.
\end{exmp}

\section*{Acknowledgments} 
This research was conducted at Warwick University as part of the author's PhD thesis and then at TU Eindhoven, where the author was supported by Jan Draisma's Vidi grant from the Netherlands Organisation for Scientific Research (NWO). The author is grateful to Franz Lehner, Diane Maclagan,  Kristian Ranestad, Jim Q. Smith, Bernd Sturmfels, the anonymous referees and the editors for helpful comments and discussions.

\end{document}